\magnification 1200
\input amssym.def
\input amssym.tex
\parindent = 40 pt
\parskip = 12 pt
\font \heading = cmbx10 at 12 true pt
 at 22 true pt
 at 19 true pt
 at 7 true pt
\def \R{{\bf R}}

\centerline{\heading An Elementary Coordinate-Dependent Local Resolution}
\centerline{\heading of Singularities and Applications}
\rm
\line{}
\line{}
\centerline{\heading Michael Greenblatt}
\line{}
\line{}

\noindent{\bf 1. Introduction}

There are many contexts in analysis and other areas of mathematics where having an explicit
and elementary resolution of singularities algorithm is helpful in understanding local properties of 
real-analytic functions, or proving theorems that depend on local properties of real-analytic
functions. In this paper, a geometric classical analysis resolution of singularities algorithm is developed. 
It is elementary in its statement and proof, heavily using explicit coordinate systems.
As one might expect, the trade-off for such a method is a weaker theorem than Hironaka's
work [H1]-[H2] or its subsequent simplications and extensions such as [BM2] [EV1] [K]
[W]. But, as will be seen, despite being entirely elementary and self-contained the theorem still can be used to prove 
various analytic results of interest.
 In addition, in [G5] (and hopefully in other future work) some of the methods
of this paper, as opposed to any specific resolution of singularities theorems, are used to give results
regarding oscillatory integrals and related objects. Although there have been other elementary approaches
to local resolution of singularities (e.g. [BM1] [Su]), the proofs here are quite different and
the algorithm has new properties tailored to classical analysis applications such as Theorem 1.1. 
This is described in more detail at the end of this section.

The author is especially interested in critical integrability
exponents, oscillatory integrals, and other objects defined through integrals whose analysis is
amenable to resolution of singularities methods. 
As illustration of our methods, two theorems are proven. First and most 
notably, a general theorem regarding the existence of critical integrability exponents is 
established. Hironaka's theory [H1] [H2] can be used to prove this result; the point here is that 
this is a new elementary classical analysis method of reasonable length that can be used to prove 
these things in short order. 
Secondly, another proof of a well-known inequality of Lojasiewicz [L] is given. Because they
are pointwise inequalities and do not depend on quantities defined through integrals, Lojasiewicz-type 
inequalities are generally easier to deal with and should be expected to follow from a reasonable
resolution of singularities procedure. In a
separate paper [G4], as another application of our methods, we give a theorem regarding asymptotic 
expansions of subvolume integrals. This result gives as a corollary the existence of asymptotic 
expansions for
oscillatory integrals (normally proved using a strong version of Hironaka's results), as well as a new
proof of the well-known result of 
Atiyah [At] and Bernstein-Gelfand [BGe] concerning the meromorphy of integrals of $f^z$ for 
real analytic $f$. 

The arguments of this paper are entirely self-contained other than using the implicit function 
theorem, the Heine-Borel Theorem, and some elementary linear algebra. No concepts from algebraic
geometry are used, not even the Weierstrass preparation theorem. The methods here can be viewed as 
generalizing the two-dimensional algorithm of [G1], as well as the papers [PS] and [V]. As one might 
expect however, the two-dimensional argument is substantially simpler. On a more
technical level, some considerations from [BM1] and [BM3] were useful in generalizing to 
$n$ dimensions. In turn, [BM1] and [BM3] 
are very much related to Hironaka's monumental work [H1] and [H2]. It should be pointed out 
that there has also been much important recent work in this area on the algebraic side. For example 
Encinas and Villamayor [EV1]-[EV2], Kollar [K], and Wlodarczyk [W] have recently given general 
resolution of singularities theorems in a more abstract setting.

To motivate our theorem concerning the existence of critical integrability 
exponents, suppose $f(x)$ is a real-analytic function on a neighborhood of the origin 
such that $f(0) = 0$. For a bounded open set 
$U$ containing $0$ and for $\epsilon > 0$ define the quantity $I_U^{\epsilon}f$ by
$$I_U^{\epsilon}f = \int_U |f|^{-\epsilon} \eqno (1.1)$$ 
For any given $U$, if $I_U^{\epsilon}f < \infty$ one automatically has that $I_U^{\epsilon'}f
< \infty$ for $\epsilon' < \epsilon$; this follows for example from Holder's inequality.
On the other hand, if $\epsilon$ is large enough then $I_U^{\epsilon}f = 
\infty$. Thus there is some critical $\epsilon$, call it $\epsilon_0$, such that 
$I_U^{\epsilon}f < \infty$ for $\epsilon < \epsilon_0$ and $I_U^{\epsilon}f = \infty$ 
for $\epsilon > \epsilon_0$. Such an $\epsilon_0$ is called a ``critical integrability exponent'' by 
many analysts, and is related to what is called a ``log-canonical threshold'' by
many people working in algebraic geometry. Besides being of intrinsic interest, such quantities
comes up frequently in geometric analysis (see [PSt], [T] for example) and it 
has long been understood resolution of singularities plays a major role in their analysis.
Also, one might ask if at the critical $\epsilon_0$ do we necessarily have
$I_U^{\epsilon}f = \infty$ for $\epsilon = \epsilon_0$. 
The answer to this question is yes; it is a consequence of the case $m=1$ and $K = \{0\}$ of the 
following theorem:

\noindent {\bf Theorem 1.1:} Suppose $f_1(x),...,f_m(x)$ are real-analytic functions defined on a
neighborhood $V$ of a compact subset $K$ of $\R^n$. There is a neighborhood $V'$ of $K$ with
$K \subset V' \subset V$ and finitely
many inequalities $\sum_{j=1}^m s_{ij} \epsilon_j < t_i$ such that if $O$ is an open set with $K
\subset O \subset V'$, then $\int_O \prod_{l=1}^m |f_l|^{-\epsilon_l} <  \infty$ if and only if
$\sum_{j=1}^m s_{ij} \epsilon_j < t_i$ for each $i$. Here the $s_{ij}$ are nonnegative
rational numbers and the $t_i$ are positive rational numbers, all independent of $O$.

Note that Theorem 1.1 is trivial if one of the functions $f_l(x)$ is the zero 
function, or if all $f_l(x)$ are nonvanishing on $K$. So the relevant situation 
is when at least one of the functions has a zero on $K$ but none of the functions is the zero
function. Also, the $\epsilon_j$ in Theorem 1.1 do not all have to be positive. 

\noindent The inequality of Lojasiewicz that we will prove is as follows ([L]):

\noindent {\bf Lojasiewicz Inequality:} Suppose $K$ is a compact set, and $f_1$ and $f_2$ are 
real-analytic functions on an open set $V$ containing $K$ such that
$\{x \in V: f_2(x) = 0\} \subset \{x \in V: f_1(x) = 0\}$. Then there is an open set $V'$ with
$K \subset V' \subset V$ and constants $\mu, C > 0$ such that $|f_2| \geq C|f_1|^{\mu}$ on $V'$. 

In the Main Theorem, an arbitrary bump function on a neighborhood of the origin will be written as the 
sum of finitely many functions. Each of these functions, after the coordinate changes given by the 
Main Theorem, becomes a quasibump function as defined below. As explained at the end of this 
section, quasibump functions are amenable to integrations by parts such as when dealing with 
oscillatory integrals. In the sequel [G4] to this paper, it is shown that
rather than having a quasibump function in the blown-up coordinates, one can just have a smooth function
times the characteristic function of the product of $n$ intervals. 
However, showing this requires a fair amount
of additional argument and so we refer to that preprint for details.

\noindent {\bf Definition:} Let $E = \{x : x_i > 0$ for all $i\}$ and let $\bar E$ denote its
closure. If $h(x)$ is 
a bounded, nonnegative, compactly supported function on $E$, we say $h(x)$ 
is a {\it quasibump function} if $h(x)$ is of the following form:
$$h(x) = a(x) \prod_{i=1}^j b_i (c_i(x) {p_i(x) \over q_i(x)}) \eqno (1.2)$$
Here $p_i(x), q_i(x)$ are monomials, $a(x) \in C^{\infty}\bar E$, the $c_i(x)$ are  
nonvanishing real-analytic functions defined
on a neighborhood of $supp(h)$, and $b_i(x)$ are nonnegative
functions in $C^{\infty}(\R)$ such that there are  $c_1 > c_0 > 0$ with each $b_i(x)
= 1$ for $x < c_0$ and $b_i(x) = 0$ for $x > c_1$. 

\noindent We now define the two key types of coordinate changes used in this paper.

\noindent {\bf Definition:} We call a function $m:A \subset \R^n \rightarrow \R^n$ an 
{\it invertible monomial map} if there are nonnegative
integers $\{\alpha_{ij}\}_{i,j=1}^n$ such that the matrix $(\alpha_{ij})$ is invertible 
and $m(x) = (m_1(x),...,m_n(x))$ where $m_i(x) = x_1^{\alpha_{i1}}....x_n^{\alpha_{in}}$.
The matrix $(\alpha_{ij})$ being invertible ensures that $h$ is a bijection on $\{x : x_l > 0
\hbox { for all } l\}$. 

\noindent {\bf Definition:} We say that a function $g: A \subset \R^n \rightarrow \R^n$ a {\it quasi-translation}
if there is a real analytic function $r(x)$ of $n-1$ variables such that 
$g(x) = (g_1(x),...,g_n(x))$, where for some $j$ we have $g_j(x) = x_j - r(x_1,...
x_{j-1},x_{j+1},...,x_n)$ and where $g_i(x) = x_i$ for all $i \neq j$.
 In other words $g(x)$ is a translation in the $x_j$ variable when the others are fixed.

In this paper, the functions we will need to resolve the zero set of a function will all be
reflections, translations, invertible monomial maps, and quasi-translations. 
The invertible monomial maps here serve the traditional purpose of blow-ups in a resolution
of singularities process such as [H1]-[H2]. However, the monomial
maps appearing in this paper will not necessarily be blow-ups or finite compositions of 
blow-ups; the type of geometric arguments used here require a broader collection
of coordinate changes. The purpose of the quasi-translations will be to convert surfaces
defined by the implicit function theorem into hyperplanes. 

\noindent We now come to the main theorem of this paper, giving Theorem 1.1 as well as the 
Lojasiewicz inequality mentioned earlier. As mentioned before, Hironaka's results and their more
recent simplifications/extensions give more general resolution of singularities theorems; the goal here is to
develop a new analytic technique (not just theorems), explicit and entirely elementary, that is applicable to 
various situations in classical analysis such as those given here and in the sequels [G4] and [G5] to 
this paper.

\noindent {\bf Main Theorem:} Let $f(x)$ be a real-analytic function defined in a neighborhood
of the origin in $\R^n$. Then there is a neighborhood $U$ of the origin such that if $\phi(x) 
\in C_c^{\infty}(U)$ is nonnegative with $\phi(0) > 0$, then 
$\phi(x)$ can be written (up to a set of measure zero) as a finite sum $\sum_i \phi_i(x)$ of
nonnegative functions such that for all $i$, $0 \in supp(\phi_i)$ and $supp(\phi_i)$ is a subset of 
one of the $2^n$ closed quadrants defined by the hyperplanes $\{x_m = 0\}$. The following 
properties hold.

\noindent (1) For each $i$ there are bounded open sets $D_i^0$,...,$D_i^{k_i}$, and maps 
$g_i^1$,..., $g_i^{k_i}$, each a reflection, translation, invertible monomial map, or quasi-translation, 
such that $D_i^0 = \{x: \phi_i(x) > 0\}$ and such that each $g_i^j$ is a real-analytic 
diffeomorphism from $D_i^j$ to $D_i^{j-1}$. The function $g_i^j$ extends to 
a neighborhood $N_i^j$ of the closure $\bar D_i^j$ with $g_i^j(N_i^j) \subset N_i^{j-1}$ for $j > 1$ and
$g_i^1(N_i^1) \subset U$. \parskip=3pt

\noindent (2) Let $E = \{x : x_i > 0$ for all $i\}$ and $\Psi_i = g_i^1 \circ .... \circ g_i^{k_i}$.
Then $D_i^{k_i} \subset E $, and there is a quasibump function
$\Phi_i$ such that $ \chi_{D_i^{k_i}}(x)(\phi_i \circ \Psi_i(x)) =  \Phi_i(x)$.

\noindent (3) $0 \in N_i^{k_i}$ with $\Psi_i(0) = 0$.

\noindent (4) On $N_i^{k_i}$, the functions $f \circ \Psi_i$, $det(\Psi_i)$, and each $j$th component 
function $(\Psi_i)_j$ is of the form $c(x)m(x)$, where $m(x)$ is a  monomial and $c(x)$ is
nonvanishing.
\parskip=12pt

To be clear, in (4) above, $det(\Psi_i)$ refers to the Jacobian determinant of $\Psi_i$. 
It is often useful to resolve several functions simultaneously, and in conjunction with
Lemma 2.3 the Main Theorem immediately gives the following corollary:

\noindent {\bf Corollary to the Main Theorem:} Suppose $\{f_l(x)\}_{l=1}^m$ are real-analytic functions
defined on a neighborhood of the origin. Then there is a neighborhood 
$U$ of the origin on which
each $f_l(x)$ satisfies the conclusions of the Main Theorem, such that for any $\phi(x)$
one can use the same decomposition $\phi = \sum_i \phi_i$ and the same coordinate changes 
$g_i^j$ to resolve each $f_l(x)$. 

To give an idea of the methods that will be used in this paper, as well as some of its 
antecedents, we turn our attention to another subject in which explicit methods of resolving 
singularities have proved useful, the study of oscillatory integral operators. Consider the 
operator on $L^2(\R)$ given by
$$T_{\lambda}f(x) = \int_{\R}e^{i \lambda S(x,y)}\chi(x,y) f(y)\,dy \eqno (1.3)$$
Here $\lambda$ denotes a parameter and $\chi(x,y)$ is a cut-off function supported near the 
origin with $\chi(0,0) \neq 0$. It is natural to ask what is the 
supremum of the $\delta$ for which $T$ has $L^2$ smoothing of order $|\lambda|^{-\delta}$. In other 
words, we seek the supremum of the $\delta$ for which there is a constant $C$ with
$$||T_{\lambda}f||_{L^2} \leq C|\lambda|^{-\delta}||f||_{L^2} \eqno (1.4)$$
It turns out that the second derivative $\partial^2 S \over \partial x\partial y$ plays a key
role. In fact, in order for there to be any $\delta$ at all for which $(1.4)$ holds there must be 
some multiindex $(\alpha,\beta)$ for which $\partial_x^{\alpha}\partial_y^{\beta}
{\partial^2 S \over \partial x\partial y}(0,0) \neq 0$.
In the case of real-analytic $S(x,y)$, in [PS] Phong and Stein found a general expression 
for the best $\delta$ for which 
$(1.4)$ holds. Their formula was in terms of the Newton polygon of $\partial^2 S \over \partial x\partial y
$ at $(0,0)$, and in their proof they divided a neighborhood of the origin into 
``curved sectors'' in
two dimensions that arise from applying the Weierstrass Preparation Theorem to ${\partial^2 S \over 
\partial x\partial y}$. Thus these sectors derive from the
resolution of singularities of ${\partial^2 S \over \partial x\partial y}$. One reason it is
useful
to divide a neighborhood in this way is that one may do different coordinate changes on each
curved sector to get a function into a desirable form. In [PSSt] related
concepts are used in their study of oscillatory integral analogues of $(1.1)$, again in two dimensions. 
With the goal of developing methods that were hoped to be generalizable to any number of 
dimensions where Puiseux-type expansions and preparation theorems are hard to find, the author 
devised effective resolution of singularities algorithms in two dimensions and in [G1] reproved the 
main result of [PS], again using appropriate curved sectors, this time coming from the resolution of 
singularities algorithm. Generalizations to general $C^{\infty}$ functions are proven in
[G2] (see also [R]), and applications to quantities of the form $(1.1)$ are in [G3].

Thus inspired, for the purpose of analyzing critical integrability exponents and oscillatory
integral operators, instead of 
requiring a single sequence of coordinate changes to work on a whole neighborhood of the 
origin, it makes sense to try to take a nonnegative bump function $\phi$ equal to 1
on a neighborhood of the origin,
and write it as $\phi = \sum_i \phi_i$. We want there to be
a sequence of canonical coordinate changes on the support of $\phi_i$ whose composition $\Psi_i$ 
takes $0$ to $0$ and converts $f$ into a monomial times a nonvanishing function. One
also wants that $\Psi_i$ is one-to-one on $\Psi_i^{-1} \{x: \phi_i(x) > 0\}$ and that the Jacobian 
of $\Psi_i$ is comparable to a monomial. The critical integrability 
exponent of $f(x) \phi_i(x) $ is then given by a rational number: If $f$ times the Jacobian
of $\Psi_i$ is written as $c(x)x_1^{\alpha}....x_n^{\alpha_n}$ in the final coordinates,
$|c(x)|$ bounded away from 0, then assuming we are integrating over a bounded set containing 
some box $(0,\eta)^n$ on which $\phi_i \circ \Psi_i > 0$, the  exponent is just ${1 \over \max_i\alpha_i}$.
Consequently the critical integrability exponent for $f \phi = \sum_i f \phi_i$ is the minimum of 
these numbers over all $i$. 

An argument given in section 5 allows this idea to be extendable 
to the setting of Theorem 1.1. It should be pointed out
that the idea of partitioning a neighborhood of the origin into such curved sectors according 
to the singularities of $f$ has before been also been used in the study of objects such as 
$(1.1)$ in some relatively nondegenerate settings; it comes up when one uses the Newton 
polyhedron of $f$ to understand the growth of $|f|$ near the origin. (See [V] for example.)

Some general heuristics behind the proof of the Main Theorem are as follows. We will proceed by 
induction on the 
dimension $n$, and for a given $n$ we induct on the order $k$ of the zero of $f$ at the origin. 
In section 2, we will
prove some technical lemmas and then reduce consideration to a class of functions amenable to
the methods of this paper. Then in section 3, we will subdivide
the cube $(-\delta,\delta)^n$ into finitely many pieces. On each piece, the terms 
of $f(x)$'s Taylor expansion corresponding to a particular face or vertex of
the Newton polyhedron of $f$ ``dominate''. Verifying that this subdivision has the 
requisite properties will take up most of section 3. In section 4 we do some further 
subdivisions of
these pieces so that, after a few allowable transformations, each subpiece becomes a set $D$
such that $(0,\rho)^n \times V \subset D \subset (0,\rho')^n \times V$ for some $\rho' > \rho > 0$
and some open set $V$ not intersecting any of the hyperplanes $\{x_m = 0\}$. Furthermore, 
under the composition of these transformations,
$f(x)$ transforms into a function of the form $m(x)h(x)$, where $m(x)$ is a monomial and where 
$h(x)$ typically has a zero of order $l < k$ at the origin. For a $\phi \in
C_c(-\delta,\delta)^n$, one writes $\phi = \sum_i \alpha_i$, where the support of $\alpha_i(x)$ 
converts under these transformations into a function 
approximately supported on the associated subpiece $D$. 

The fact that $l < k$ allows one to iterate the above idea, further subdividing the subpieces and 
further decomposing the $\alpha_i(x)$ until one
finally has the $\phi_i(x)$ satisfying the conclusions of the Main Theorem. Lemma 2.2 is an 
important technical lemma that facilitates this induction step.

We now draw attention to some properties of the Main Theorem, tailored towards applications, that are 
not present in other elementary resolution of singularities methods. First, part (4) of 
the Main Theorem gives that the determinant of the composition $\Psi_i$ of the coordinate changes is
comparable to a monomial, a key fact used in section 5 when proving Theorem 1.1. Also, the form
$(1.2)$ of the function $\Phi_i$
in $(2)$ of the Main Theorem is appropriate for oscillatory integrals and related objects; if 
$\Phi_i$
were too irregular then integrations by parts in the blown-up domain can cause unnecessarily 
large factors to show up. As it is here, an $x_k$ derivative landing on $\Phi_i$ only
gives an additional factor of ${C \over x_k}$, appropriate for most purposes. On a related note, the
fact that the $D_i^{k_i}$ are disjoint and the $\Psi_i$ are one-to-one (and not some $m$ to one) on 
$D_i^{k_i}$ ensures that 
in pulling back integrals to blown-up coordinates one does not to do any unnatural subdivisions of 
$D_i^{k_i}$ which may again lead to problems doing integrations by parts or analyzing critical 
integrability exponents. It is also worth pointing out that by part 4) of the Main Theorem, 
the $f \circ \Phi_i$ and $det(\Phi_i)$ are comparable to a monomial on a neighborhood of $\bar D_i^{k_i}$.
This further helps avoid problems at the boundaries of the $D_i^{k_i}$ when one is performing
integrations.

There are also substantial differences between the proof of the Main Theorem and other resolution of
singularities theorems, including the elementary algorithms of [BM1] and [Su]. The Newton polyhedron-based coordinate-dependent subdivisions described above
exemplify this. Like in all papers in this subject, we do have an invariant that decreases under
each iteration of the algorithm. The invariant here
is simply the order of vanishing of the function being resolved, and not a more elaborate 
invariant as in [BM1]-[BM3] [EV1]-[EV2] [H1]-[H2]. In addition, we do not have to 
consider the maximum stratum of an invariant as in [BM1]-[BM3]. As indicated above, if 
we are at a stage
of the resolution process where a function being resolved has a zero of some order $k$, one does
some subdivisions and coordinate changes and then one factors out a monomial. If the resulting function
still has a zero of order $k$, one does a quasi-translation coming from the implicit function theorem
applied to a $k-1$th derivative of the new function. In the language of resolution of singularities,
this quasi-translation takes a ``hypersurface
of maximal contact'' containing the $k$th order zeroes of the function, determined by the coordinate
system we are working in, and translates it to the
hyperplane $x_n = 0$. The zeroes of all further functions in the resolution 
process will be of order at most $k-1$. One does not have to keep
track of any history of the resolution process as in [BM1]-[BM3] [EV1]-[EV2] [H1]-[H2].
In fact, one does not have to consider exceptional hypersurfaces in any form in this paper.

\noindent {\bf 2: Beginning the proof of the Main Theorem: a localization lemma; preparation of 
the function}

We start with a relatively easy lemma saying that if the product of several functions on a 
cube is comparable to a monomial, so is each of the individual functions:

\noindent {\bf Lemma 2.1:} Suppose $N$ is an open cube with $0 \in N$ and 
$\{a_i(x_1,...,x_q)\}_{i=1}^p$ are real-analytic
functions defined on its closure $\bar N$ such that the following holds:
$$\prod_{i=1}^p a_i(x_1,...,x_q) = m(x_1,...x_q)s(x_1,...,x_q) \eqno (2.1)$$ 
Here $m(x_1,...,x_q)$ is a monomial and $s(x) \neq 0$ on $\bar N$.
Then each $a_i(x_1,...,x_q)$ can be written as $m_i(x_1,...x_q)s_i(x_1,...,x_q)$, where each 
$m_i(x)$ is a monomial and $s_i(x) \neq 0$ on $\bar N$.

\noindent {\bf Proof:} Write $m(x_1,...x_q) = cx_1^{\alpha_1}...x_n^{\alpha_n}$. We induct on
$s = \sum_j \alpha_j$. If $s = 0$ the result is trivial, so assume that $s > 0$. Let $j_0$ be
an index such that $\alpha_{j_0} > 0$. then on the hyperplane $\{x: x_{j_0} = 0\}$, 
$\prod_{i=1}^p a_i(x_1,...,x_q) = 0$. As a result, at least one of the functions $a_i(x_1,...,x_q)$,
call it $a_{i_0}(x_1,...,x_q)$, must also be the zero function on this hyperplane. So we can write
$a_{i_0}(x_1,...,x_q) = x_{j_0}a_{i_0}'(x_1,...,x_q)$ for some real-analytic function 
$a_{i_0}'(x_1,...,x_q)$. We then have
$$a_{i_0}'(x_1,...,x_q)\prod_{i\neq i_0}a_i(x_1,...,x_q) = cx_1^{\alpha_1}...
x_{j_0}^{\alpha_{j_0} - 1}...x_n^{\alpha_n}$$
The result then follows from the induction hypothesis, and we are done.

\noindent The following important ``localization lemma'' is used in the inductive step.

\noindent {\bf Lemma 2.2:} Suppose $F(x)$ is a real-analytic function defined in a neighborhood
of the origin and $\beta(x)$ is nonnegative with $\beta(x) = 1$ on 
a neighborhood of the origin.
Suppose that we can write $\beta(x) = \sum_{i=1}^r \beta_i(x)$, where for each $i$, $\beta_i(x) \geq 0$, 
$0 \in supp(\beta_i)$, and $supp(\beta_i)$ is a subset of one of the closed quadrants defined by the hyperplanes 
$\{x_m = 0\}$.
Suppose further that to each $\beta_i(x)$ there is a bounded open set $D_i \subset E = \{x: x_i > 0$ 
for all $i\}$ and a real-analytic diffeomorphism $\zeta_i(x)$ from $D_i$ to $\{x: \beta_i(x) > 0\}$ 
such that the following hold, where as usual $\bar D_i$ denotes the closure of $D_i$.

\noindent 1) $\zeta_i(x)$ extends to a continuous function on an open set $N_i$ containing 
$\bar D_i$.\parskip=3pt

\noindent 2) $\zeta_i(x)$ is a composition $g_i^1 \circ ... \circ g_i^{k_i}$ of translations, 
reflections, quasi-translations, and invertible monomial maps satisfying 1) and 4) of the Main 
Theorem with $\phi_i(x) = \beta_i(x)$ and $N_i^{k_i} = N_i$. In particular, on $N_i$, 
$det(\zeta_i)$ and each $j$th 
coordinate function $(\zeta_i)_j(x)$ is of the form $c(x)m(x)$, where $m(x)$ is a 
monomial and $c(x)$ is nonvanishing.  

\noindent 3) There is a quasibump function $\Phi_i(x)$ such that $\chi_{D_i}(x)(\beta_i \circ \zeta_i(x) )
= \Phi_i(x)$.

\noindent 4) For each $w \in \bar D_i \cap \zeta_i^{-1}(0)$, there is a quasi-translation $q_w$ in 
the $x_n$ variable with $q_w(w) = w$ such that $F \circ \zeta_i \circ q_w (x + w)$
satisfies the conclusions of the Main Theorem in some neighborhood of $x = 0$. If $w_n = 0$, 
assume $q_w$ is the identity map. 
\parskip = 12pt

\noindent Then the Main Theorem holds for $F(x)$ on some neighborhood of the origin. 

\noindent {\bf Proof:} Let $K = \bar D_i \cap \zeta_i^{-1}(0)$. Suppose $w \in K$. 
By assumption, we may let $U_{i,w}$ be an open set on which $F \circ \zeta_i \circ q_w(x + w)$ 
satisfies the conclusions of the Main Theorem. Let $\eta^{i,w}(x) \in C_c^{\infty}
(U_{i,w})$ be a bump function which is equal to 1 on some neighborhood $V_{i,w}$ of $0$. 
Then $w \in q_{w}(V_{i,w} + w)$. By compactness of $K$ we may
let $\{w_{ij}\}$ a finite set of points such that $K \subset \cup_j \,q_{w_{ij}}(V_{i,w_{ij}} + 
w_{ij})$. 

Let $\eta^{i,w} = \sum_l \eta^{i,w}_l$ be the decomposition from the Main Theorem. We 
use the $\eta^{i,w_{ij}}_l$ to give a sort of partition of unity on $\cup_j \,q_{w_{ij}}(V_{i,w_{ij}}
+ w_{ij})$; namely we let $\gamma^{ijl}(x) = {\eta^{i,w_{ij}}_l(q_{w_{ij}}^{-1}x 
- w_{ij})  \over \sum_{j,l} \eta^{i,w_{ij}}_l(q_{w_{ij}}^{-1}x - w_{ij})}= 
{\eta^{i,w_{ij}}_l(q_{w_{ij}}^{-1}x - w_{ij})  \over \sum_j \eta^{i,w_{ij}}(q_{w_{ij}}^{-1}x - 
w_{ij})}$. 
If $\Psi_{ijl}^2$ is the composition of the coordinate changes from the Main Theorem, then by 2) of 
the Main Theorem for a $\Phi_{ijl}(x) \in C^{\infty}(E)$ we have
$$\eta^{i,w_{ij}}_l \circ \Psi_{ijl}^2(x) =  \Phi_{ijl}(x) \eqno (2.2a)$$ 
Adjusting coordinates, this implies that we can similarly write 
$$\gamma^{i,w_{ij}}_l \circ q_{w_{ij}} (\Psi_{ijl}^2(x) + w_{ij})  =  \Phi_{ijl}'(x) 
\eqno (2.2b)$$
As before $\Phi_{ijl}'$ is in $C^{\infty}(E)$. Another useful 
observation is the following. By the Main Theorem, each $\eta^{i,w}_l(x)$ is supported in one
of the $2^n$ cubes defined by the hyperplanes $x_m = 0$, so the same property
holds for each $\gamma^{ijl}(q_{w_{ij}}(x + w_{ij})) = 
{\eta^{i,w_{ij}}_l(x) \over \sum_j \eta^{i,w_{ij}}(x)}$. 
Furthermore, assuming $U_{w_{ij}}$ was chosen to be small enough, the same property holds for 
$\gamma^{ijl}(x)$. (For this we use that $q_w$ is a quasi-translation in the $x_n$ variable which is
the identity map when $w_n = 0$). As a result we have 
$$supp(\gamma^{ijl}) \subset \bar E \hbox{ or } supp(\gamma^{ijl}) \cap E = \emptyset \eqno (2.3)$$
Next, observe that for a sufficiently small $\delta$ we have that
$$\bar D_i \cap \zeta_i^{-1}[-\delta,\delta]^n \subset \cup_j \,q_{w_{ij}}(V_{i,w_{ij}} + w_{ij})
\eqno (2.4)$$
To see why $(2.4)$ holds, suppose not. Then for each $\delta > 0$ the compact set $L_{\delta}$ is 
nonempty, where
$$L_{\delta} = \bar D_i \cap \zeta_i^{-1}[-\delta,\delta]^n
\cap (\cup_j \,q_{w_{ij}}(V_{i,w_{ij}} + w_{ij}))^c \eqno (2.5)$$
Taking the intersection of $(2.5)$ over all $\delta$ would give that $L_0$ is nonempty as well,
contradicting the definition of the $V_{i,w_{ij}}$. Hence $(2.4)$ holds. Fix some $\delta > 0$
small enough that $(2.4)$ holds for each $i$, small enough that any $supp(\gamma^{ijl} \circ
\zeta_i^{-1})$ intersecting $[-\delta,\delta]^n$ contains the origin, and small enough that $\beta(x) = 1$ on 
$(-\delta,\delta)^n$. Let $\phi(x)$ be any function in $C_c^{\infty}((-\delta,\delta)^n)$; 
we will see that the Main Theorem holds for $\phi(x)$. This will imply the lemma we are proving. 

\noindent If $\phi_i$ denotes $\phi\beta_i$ we have
$$\phi = \phi\beta = \sum_i\phi\beta_i= \sum_i \phi_i$$
Furthermore, by $(2.4)$ we have
$$supp(\phi_i \circ \zeta_i)  \subset \cup_j \,q_{w_{ij}}(V_{i,w_{ij}} + w_{ij}) \eqno (2.6)$$
We may decompose a given $\phi_i(x)$ by 
$$\phi_i(x) = \sum_{j,l} \phi_i(x)\gamma^{ijl}(\zeta_i^{-1}(x)) \eqno (2.7a)$$
The sum $(2.7a)$ makes sense for the following reason. If $\phi_i(x) \neq 0$, then
$\beta_i(x) \neq 0$ and by the assumptions of this lemma, $x = \zeta_i(y)$ for a unique $y \in D_i$.
This $y$ is in $\cup_j \,q_{w_{ij}}(V_{i,w_{ij}} + w_{ij})$ by $(2.6)$, and therefore $\sum_{j,l} 
\gamma^{ijl}(y) = 1$. We conclude that the right-hand side of $(2.7a)$ makes sense and equality
$(2.7a)$ holds. We analogously decompose $\phi(x)$ via
$$\phi(x) = \sum_{i,j,l} \phi_i(x)\gamma^{ijl}(\zeta_i^{-1}(x)) \eqno (2.7b)$$
We will see that the decomposition $(2.7b)$ satisfies the conclusions of the Main Theorem, where
the maps called $g_i^j$ in the Main Theorem, denoted here by $g_{ijl}^m$, are as follows.
Write $\zeta_i = G_i^1 \circ ... \circ G_i^{k_i}$, where each $G_i^j$ is a translation, 
reflection, invertible monomial map, or quasi-translation. Let $h_{ijl}^m$ denote the analogous
transformations for $F \circ \zeta_i \circ q_{w_{ij}}(x + w_{ij})$ corresponding to $\eta^
{i,w_{ij}}_l$. 
For $1 \leq m \leq k_i$, define $g_{ijl}^m = G_i^j$. Let $g_{ijl}^{k_i + 1} = q_{w_{ij}}$, 
let $g_{ijl}^{k_i + 2}$ denote the shift
$x \rightarrow x + w_{ij}$, and for $m > k_i + 2$ let $g_{ijl}^m = h_{ijl}^{m - k_i - 2}$.

We now verify the various conclusions of the Main Theorem. First, each term in $(2.7b)$ is 
supported in one of the quadrants defined by the hyperplanes $x_m = 0$ since the same property holds
for the $\beta_i(x)$. Next, by definition each $g_{ijl}^m$
is either a translation, reflection, invertible monomial map, or quasi-translation. That each is a real-analytic
diffeomorphism on its domain, extending to a real-analytic function on a neighborhood of its closure,
follows from the corresponding property of the $G_i^m$, $h_{ijl}^{m - k_i - 2}$, $q_{w_{ij}}$, or 
shift; 
the domain of $g_{ijl}^m$ is a subset of that of the appropriate function. This gives 1) of the 
Main Theorem. As for 2), let $\Psi_{ijl}$ denote the composition of all the $g_{ijl}^m$.  
As in $(2.2a)-(2.2b)$ let $\Psi_{ijl}^2$ denote the composition of all the $h_{ijl}^m$.
Then we have
$$\Psi_{ijl}(x) = \zeta_i \circ q_{w_{ij}} (\Psi_{ijl}^2(x) + w_{ij}) \eqno (2.8)$$
Next, if we let $\phi_{ijl}(x)$ denote the term $\phi_i(x)\gamma^{ijl}(\zeta_i^{-1}(x))$ of 
$(2.7)$, observe that
$$\phi_{ijl} \circ \Psi_{ijl} (x) = \phi_i \circ \Psi_{ijl}(x)
(\gamma^{ijl} \circ q_{w_{ij}}(\Psi_{ijl}^2(x) + w_{ij}))
= \phi_i \circ \Psi_{ijl}(x)  \Phi_{ijl}'(x) \eqno (2.9)$$
The last equality follows from $(2.2b)$. Next observe that 
$$\phi_i \circ \Psi_{ijl}(x) = \phi_i \circ \zeta_i \circ q_{w_{ij}}(\Psi_{ijl}^2(x) 
+ w_{ij}) \eqno (2.10)$$
Hence 
$$ \phi_{ijl} \circ \Psi_{ijl} (x) = \phi_i \circ \zeta_i \circ q_{w_{ij}}(\Psi_{ijl}^2(x) + w_{ij})
(\gamma^{ijl} \circ q_{w_{ij}}(\Psi_{ijl}^2(x) + w_{ij})) \eqno (2.11)$$
In order for $(2.11)$ to be nonzero, $y = q_{w_{ij}}(\Psi_{ijl}^2(x) + w_{ij})$ must satisfy 
$\gamma^{ijl}(y) \neq 0$. By $(2.3)$, for a given $i$, $j$, and $l$ there are two possibilities. 
First, $q_{w_{ij}}(\Psi_{ijl}^2(x) + w_{ij})$ may never be in $\bar E$, in which case by condition 
3) of this lemma,
$\phi_i \circ \zeta_i \circ q_{w_{ij}} (\Psi_{ijl}^2(x) + w_{ij})$ is necessarily the zero function.
(Recall $\phi_i = \phi \beta_i$). But this would imply that $(2.11)$ is the zero function, 
a contradiction. So we must have the second possibility, which
is that $y = q_{w_{ij}}(\Psi_{ijl}^2(x) + w_{ij})$ is always in $E$. In this case, condition 3) says
that $\phi_i \circ \Psi_{ijl}(x) = $ $\phi_i \circ \zeta_i \circ q_{w_{ij}}(\Psi_{ijl}^2(x) + w_{ij})$ 
is a smooth function, defined on a neighborhood of the support of $\Phi_{ijl}'(x)$. Therefore, in 
view of $(2.10)$, equation $(2.9)$ 
says that $\phi_{ijl} \circ \Psi_{ijl} (x)$ is equal to $\Phi_{ijl}(x)$ for some 
$\Phi_{ijl}(x) \in C^{\infty}(E)$. We show that $\Phi_{ijl}(x)$ in fact satisfies $(1.2)$ as required 
after the proof of part 4) of this lemma below. 

Moving on now to 3), let $A_{ijl}$ denote the domain of $\phi_{ijl} \circ \Psi_{ijl}$, called 
$D_i^{k_i}$ in the statement of the Main Theorem. By part 1), the function $\Psi_{ijl}$ extends to 
an open set containing its closure $\bar A_{ijl}$. We will now show it also contains the origin. By assumption 
3) of this lemma, $\zeta_i$ extends to a 
neighborhood of $w_{ij}$, which in turn implies that $\zeta_i \circ q_{w_{ij}}$ 
also extends to a neighborhood of $w_{ij}$. Hence $\zeta_i \circ q_{w_{ij}} 
(w_{ij} + x)$ extends to a neighborhood of $0$, so $\Psi_{ijl}(x) =  
\zeta_i \circ q_{w_{ij}} (\Psi_{ijl}^2(x) + w_{ij})$ extends to a  
neighborhood of $(\Psi_{ijl}^2)^{-1}(0)$. By the Main Theorem, using assumption 4) of this lemma,
$\Psi_{ijl}^2(0) = 0$, so this neighborhood must contain the origin, and
we conclude that $\Psi_{ijl}$ extends
to a neighborhood of the origin as well. 

\noindent Furthermore, again using the Main Theorem on $\Psi_{ijl}^2$, we have that   
$$ g_{ijl}^1 \circ ... \circ g_{ijl}^{m_{final}}(0) = g_{ijl}^1 \circ ... \circ g_{ijl}^{k_i + 2}
(0) $$
This in turn is equal to 
$$ g_{ijl}^1 \circ ... \circ g_{ijl}^{k_i + 1}
(w_{ij}) = g_{ijl}^1 \circ ... \circ g_{ijl}^{k_i}(w_{ij}) = 0 \eqno (2.12)$$
The last equality of $(2.12)$ follows from the fact that each $w_{ij}$ is in $\zeta_i^{-1}(0)$.
Thus we are done proving 3) of the Main Theorem.

We next verify part 4) of the Main Theorem. Observe that by assumption 4) of this lemma,
$F \circ \Psi_{ijl}(x) = F \circ \zeta_i \circ q_{w_{ij}} (\Psi_{ijl}^2(x) + w_{ij})$  
is a monomial times a nonvanishing function on a neighborhood of the origin. Hence $F \circ 
\Psi_{ijl}$ satisfies the requirements of 4) of the Main Theorem. As for the required 
conditions on $det(\Psi_{ijl})$ and the component functions $(\Psi_{ijl})_m$, observe that 
$\Psi_{ijl} = \zeta_i \circ q_{w_{ij}} \circ s_{ij} \circ \Psi_{ijl}^2$, where $s_{ij}$ denotes the 
shift 
$x \rightarrow x + w_{ij}$. Suppose we show that each of $\zeta_i$, $q_{w_{ij}}$, $s_{ij}$, and 
$\Psi_{ijl}^2$ satify the determinant and component function conditions of part 4) of the 
Main Theorem. Then clearly their composition $\Psi_{ijl}$ will satisfy the component function
conditions, and furthermore by the chain rule, $\Psi_{ijl}$ will satisfy the determinant condition 
as well. 

Thus it suffices to 
show that the determinant and component function conditions hold for each of $\zeta_i$, $q_{w_{ij}}$, 
$s_{ij}$,
and $\Psi_{ijl}^2$. By assumption 2) of this lemma, they hold for $\zeta_i$. Since $\Psi_{ijl}^2$ is
a composition of functions arising from the Main Theorem, they hold for $\Psi_{ijl}^2$ as well.
The determinant condition holds for $s_{ij}$ since it is a shift, and for $q_{w_{ij}}$ since it is 
a quasi-translation. So it remains to show the component conditions for $q_{w_{ij}}$ and $s_{ij}$.

We start with $q_{w_{ij}}$. If $m \neq n$, then since $q_{w_{ij}}$ is a quasi-translation in the
$x_n$ variable we have $(q_{w_{ij}}(x))_m = x_m$, which is a monomial. If $m = n$, then 
if $(w_{ij})_n = 0$
by assumption $q_{w_{ij}}$ is the identity and $(q_{w_{ij}}(x))_n = x_n$, a monomial. If 
$(w_{ij})_n \neq 0$, then
if $U_{w_{ij}}$ was chosen to have diameter $\epsilon |(w_{ij})_n|$ for a small enough $\epsilon$,
then since $q_{w_{ij}}(w_{ij}) = w_{ij}$, on $U_{w_{ij}} + w_{ij}$ we have $(q_{w_{ij}}(x))_n
\sim x_n \sim (w_{ij})_n$. Hence $(q_{w_{ij}}(x))_n$ is comparable to the constant monomial 1.
Thus we have shown the component conditions for $q_{w_{ij}}$.

Moving on to $s_{ij}$, if some $m$th component $(w_{ij})_m$ of $w_{ij}$ is zero, then the 
$m$th component of $x + w_{ij}$ is exactly $x_m$ and the $m$th component of the shift is 
comparable to a monomial as required. In the case where $(w_{ij})_m \neq 0$, so long as we had 
chosen $U_{w_{ij}}$ such that the diameter of $U_{w_{ij}}$ is less than ${1 \over 2}|(w_{ij})_m|$,
then on $U_{w_{ij}}$ the component function $(x + w_{ij})_m$ satisfies  $(x + w_{ij})_m \sim 
(w_{ij})_m$ and therefore $(x + w_{ij})_m$ is comparable to the constant monomial 1. We conclude $s_{ij}$ 
satisfies the component conditions of part 4) of the Main Theorem. As a result, we have now proven 
that $\Psi_{ijl}$ satisfies the conditions required for part 4) of the Main Theorem.

Finally, we prove that $\Phi_{ijl}(x)$ satisfies $(1.2)$. For this, we first observe that since 
$\phi_{ijl}(x) = \phi_i(x) \gamma^{ijl}(\zeta_i^{-1}(x))$, it suffices to show that $\phi_i \circ 
\Psi_{ijl}(x)$ and $\gamma^{ijl} \circ \zeta_i^{-1}\circ \Psi_{ijl}(x)$ both satisfy $(1.2)$.
We start with $\phi_i \circ \Psi_{ijl}(x)$. Since $\phi_i
= \phi \beta_i$, it suffices to show $\beta_i \circ \Psi_{ijl} = \beta_i \circ \zeta_i \circ 
( q_{w_{ij}} \circ s_{ij} \circ \Psi_{ijl}^2)$  satisfies $(1.2)$. By
assumption 3) of this lemma, $\beta_i \circ \zeta_i$ is of the proper form $(1.2)$. By the proof of 
part 4) of
this lemma each component of $q_{w_{ij}} \circ s_{ij} \circ \Psi_{ijl}^2$ is of the form $c(x)m(x)$
for nonvanishing $c(x)$ and $m(x)$ a monomial. As a result, the composition $\beta_i \circ \zeta_i 
\circ (q_{w_{ij}} \circ s_{ij} \circ \Psi_{ijl}^2)(x)$ is also of the form $(1.2)$ as required.

As for $\gamma^{ijl} \circ \zeta_i^{-1}
\circ \Psi_{ijl}(x)$, observe that $\gamma^{ijl} \circ \zeta_i^{-1} \circ \Psi_{ijl}(x) = $
$(\gamma^{ijl} \circ q_{w_{ij}} \circ s_{ij}) \circ \Psi_{ijl}^2(x)$ $ = r_{ijl}(x)
(\eta_l^{i,w_{ij}} \circ \Psi_{ijl}^2(x))$ for a smooth function $r_{ijl}(x)$. The definition of 
$\eta_l^{i,w_{ij}}$ implies that by the Main Theorem the function $\eta_l^{i,w_{ij}} \circ 
\Psi_{ijl}^2(x)$ is of the form $(1.2)$.
Multiplying by a smooth factor does not change this, so $\gamma^{ijl} \circ \zeta_i^{-1}
\circ \Psi_{ijl}(x)$ also satisfies $(1.2)$. Hence we conclude that $\Phi_{ijl}(x)$
satsifies $(1.2)$ as well. This completes the proof of Lemma 2.2.

The following lemma is quite useful, for example in proving the corollary to the Main Theorem.

\noindent {\bf Lemma 2.3:} 
Suppose $\{f_l(x)\}_{l=1}^m$ are real-analytic functions such that $F(x) = \prod_{l=1}^m f_l(x)$ 
satisfies the conclusions of the Main Theorem on a neighborhood $U$ of the origin. Then
each $f_l(x)$ also satisfies the conclusions of the Main Theorem on $U$, such that for any $\phi(x)$
one can use the same decomposition $\phi = \sum_i \phi_i$ and the same coordinate changes 
$g_i^j$ to resolve each $f_l(x)$. 

\noindent {\bf Proof:} Let $\beta(x) \in C_c^{\infty}(U)$ be a nonnegative function with 
$\beta(x) = 1$ on a neighborhood of the origin. Let 
$\beta(x) = \sum_i \beta_i(x)$ denote the decomposition given by the Main Theorem applied to $F(x)$. Let
$g_i^j$ denote the associated coordinate changes, and let  $\zeta_i$ denote the 
composition $g_i^1 \circ ... \circ g_i^{k_i}$. By 4) of the Main Theorem, 
if $w \in supp(\beta_i \circ \zeta_i)$, the function $F \circ \zeta_i$ is a monomial 
times a nonvanishing real-analytic function on a neighborhood of $w$. Hence by Lemma 
2.1, each $f_i \circ \zeta_i$ is also a monomial times a nonvanishing real-analytic function
on some neighborhood $U_w$ of $w$. 

Shrinking $U_w$ if necessary, we assume that the diameter of 
$U_w$ is less than ${1 \over 2}\{ min_j |w_j|: w_j \neq 0 \}$. Then on $U_w - w$, the function 
$f_l \circ \zeta_i(x + w)$ is also a monomial times a nonvanishing real-analytic function. 
The reason this is true is as follows: for any $p$ for which $w_p = 0$, the $p$th component
$(x + w)_p$ is $x_p$. For any $p$ for which $w_p \neq 0$, $(x + w)_p \sim w_p$ since the 
diameter of $U_w$ is less than ${1 \over 2}\{ min_j |w_j|: w_j \neq 0 \}$. Hence under the 
coordinate change $x \rightarrow x + w$, a monomial on $U_w$ turns into a monomial times a 
nonvanishing function on $U_w - w$. We 
conclude that on $U_w - w$, each function $f_l \circ \zeta_i(x + w)$ is also a monomial times 
a nonvanishing function. We can now invoke Lemma 2.2 with $q_w$ the identity map; each 
$f_l \circ \zeta_i(x + w)$ automatically
satisfies the conlusions of the Main Theorem: one needs one coordinate change, the identity map.
Hence by Lemma 2.2 each $f_l$ also satisfies the conditions of the Main Theorem on a 
neighborhood of the origin. 

For any $\phi(x)$ the same coordinate changes work for each $f_l$; in the proof of 
Lemma 2.2 the coordinate changes are explicitly given by first the $g_i^j$,  
then the shift $x \rightarrow x + w$, then the coordinate changes on $f_l \circ 
\zeta_i(x + w)$ (only the identity map here). Furthermore in the proof of Lemma 2.2, the 
decomposition $\phi = \sum_i \phi_i$ is given in terms of the $\beta_i$, the $\zeta_i$, and 
the decompositions of bump functions induced by the coordinate changes on 
$f_l \circ \zeta_i(x + w)$, which, being the identity map are independent of $l$. This completes 
the proof of Lemma 2.3.

We now commence the proof of the Main Theorem. We prove it by induction on the dimension $n$. The 
case $n = 1$ is easy, so we assume we know the result for all dimensions less than $n$ and are
now proving it for $n$. 

The idea now is to reduce consideration to a canonical form amenable to the analysis 
of sections 3 and 4. Namely, we reduce to the case where ${\partial^{\kappa} f \over \partial 
x_n^{\kappa}}(0) \neq 0$ for some $\kappa$. We do this as follows. Let $f(x)$ be an arbitrary real-analytic 
function with a zero of some order $\kappa$ at the origin. Let $L$ be a linear map such that each 
component function $L_i$ has a nonzero $x_n$ coefficient, and such that ${\partial^{\kappa} (f 
\circ L) \over \partial x_n^{\kappa}}(0) \neq 0$. Suppose we know the Main Theorem for functions 
that vanish at the origin to finite order in the $x_n$ direction. We may apply Lemma 2.3 to $f 
\circ L$ and the functions $L_i$. Let $U$ be the associated neighborhood of the origin. I claim 
that $f$ then satisfies the conclusions of the Main Theorem on $L(U)$; if $\phi
\in C_c(L(U))$ then we may decompose $\phi \circ L = \sum \phi_i \circ L$ according to Lemma 
2.3 and let $g_i^j$ be the associated maps. For $\phi$ itself, the decomposition is then given by 
$\phi = \sum_i \phi_i$ and the associated coordinate changes are given by first some linear
quasi-translations whose composition is $L$ and then the 
sequence of $g_i^j$; this satisfies the conditions of the Main Theorem. 

Most of the conclusions of the Main Theorem follow pretty much directly from the corresponding
conclusions on $f \circ L$ and the $L_i$. A couple of things are worth pointing out. First, since
the domain $D_i$ of $L \circ g_i^1 \circ ... \circ g_i^{k_i}$ is a subset of $\{x: x_i > 0$ for 
all $i\}$ and since each component of $L \circ g_i^1 \circ ... \circ g_i^{k_i}$ is comparable to a 
monomial, we have that the image $D_i$ under
$L \circ g_i^1 \circ ... \circ g_i^{k_i}$ is a subset of one of the quadrants defined by the
hyperplanes $\{x_k = 0\}$, as required in the first paragraph of the Main Theorem. Another thing 
worth mentioning is the verification of part 4) of the Main Theorem. By assumption, $f \circ 
(L \circ g_i^1 \circ ... \circ g_i^{k_i})$ $= (f \circ L) \circ
(g_i^1 \circ ... \circ g_i^{k_i})$, is comparable to a monomial. The 
determinant of $L \circ g_i^1 \circ ... \circ g_i^{k_i}$ is a constant times the determinant
of $ g_i^1 \circ ... \circ g_i^{k_i}$, which is comparable to a monomial also. Lastly, each
$j$th component of $L \circ g_i^1 \circ ... \circ g_i^{k_i}$ is also assumed to be comparable to a 
monomial, and we conclude part 4) of the Main Theorem holds.

We conclude that we may restrict our attention to functions satisfying ${\partial^{\kappa} f \over 
\partial x_n^{\kappa}}(0) \neq 0$ for some $\kappa$. We will prove the Main Theorem for a given $n$
by induction on $\kappa$. If $\kappa = 0$ there is nothing to prove since the function is already 
comparable to a (constant) monomial, so we assume we have the Main Theorem for $\kappa - 1$ and 
are  seeking to prove it for $\kappa$. 

Next, we will further simplify the class of functions we need to consider. In fact, by an 
appropriate application of Lemma 
2.2, we will see that the inductive step of the Main Theorem follows from the following:

\noindent {\bf Theorem 2.4:} Suppose $g(x)$ is real-analytic on a neighborhood of the origin
and has Taylor expansion of the form
$$g(x) = g_{\kappa}(x_1,...x_n)x_n^{\kappa} + \sum_{l=0}^{\kappa - 1}m_l(x_1,...,x_{n-1})
s_l(x_1,...,x_{n-1})x_n^l \eqno (2.13)$$
Here $m_l(x_1,...,x_{n-1})$ is either the zero function or a nonconstant monomial, 
$g_{\kappa}(0,...,0) \neq 0$, and each $s_l(0) \neq 0$. Suppose the Main Theorem is known in dimensions
less than $n$, and in $n$ dimensions for functions that vanish to order less than $\kappa$ in the
$x_n$ direction. Then there is a cube $(-\eta,\eta)^n$ such that $g(x)$ satisfies the conclusions 
of the Main Theorem on $(-\eta,\eta)^n$.

\noindent {\bf Proof of Reduction to Theorem 2.4:}

\noindent Let $f(x)$ be a real-analytic function satisfying ${\partial^{\kappa} f \over 
\partial x_n^{\kappa}}(0) \neq 0$. We Taylor expand $f(x)$ about the origin as 
$$f(x_1,...,x_n) = f_{\kappa}(x_1,...,x_n)x_n^{\kappa} + \sum_{l < \kappa} f_l(x_1,...,x_{n-1})x_n^l 
\eqno (2.14)$$
Here $f_{\kappa}(0) \neq 0$. We also assume $f_l(0) = 0$ for $l < \kappa$; otherwise we could just
invoke the induction hypothesis for an $l < k$. We now use the induction 
hypothesis in dimension $n-1$ to simplify the form of the coefficient 
functions $f_l(x_1,...,x_{n-1})$ for $l < \kappa$. Namely, we apply Lemma 2.3 to the (nonzero) functions 
$f_l(x_1,...,x_{n-1})$ for $l < \kappa$. We get an open set $U^0 \subset \R^{n-1}$ containing
the origin satsifying the conclusions of the Main Theorem. Let $\phi \in C_c^{\infty}(U^0)$ be a bump
function which is equal to 1 on some cube $[-\delta_0,\delta_0]^{n-1}$, and let $\phi = 
\sum_i \phi_i$ be the decomposition coming from the Main Theorem. Let $g_i^j$ be the 
corresponding coordinate changes, and $D_i$ be the corresponding domains. 

Let $\beta \in C_c^{\infty}(-\delta_0,\delta_0)^n$ be a nonnegative function equal to 1
on a neighborhood of the origin. Let $\beta_i = \beta \phi_i$. The decomposition
$\beta = \sum \beta_i$, after a slight modification, will allow us to apply Lemma 2.2 and reduce
things to proving Theorem 2.4.
Define $\bar{g}_i^j(x_1,...,x_n) = (g_i^j(x_1,...,x_{n-1}),x_n)$. Then by the Main Theorem,
$\bar{g}_i^j$ is a real-analytic diffeomorphism from $\{x: \phi_i \circ 
\bar{g}_i^1 \circ ... \circ \bar{g}_i^j(x) > 0\}$ to 
$\{x: \phi_i \circ \bar{g}_i^1 \circ ... \circ \bar{g}_i^{j-1}(x) > 0\}$, and if $\zeta_i$ denotes
$\bar{g}_i^1 \circ ... \circ \bar{g}_i^{k_i}$, then $\phi_i \circ \zeta_i(x)$ is of the form
$\Phi_i(x) \chi_{E'}(x)$,  where $E'$ denotes $\{x: x_i > 0$ for $i < n\}$ and where 
$\Phi_i$ is a quasibump function in the first $n-1$ variables. Consequently, since $\beta_i = \beta 
\phi_i$, $\beta_i \circ \zeta_i$ is also of this form. 

In view of the statement of Lemma 2.2, we would like to replace $E'$ by
$E = \{x: x_i > 0$ for all $i\}$. So we write each $\beta_i = \beta_i^+ + \beta_i^-$, where
$\beta_i^+ \circ \zeta_i = \Phi_i \chi_{E}$ and $\beta_i^- \circ \zeta_i = \Phi_i \chi_{r(E)}$, 
where $r(E)$ denotes the reflection of $E$ about the hyperplane $x_n = 0$. 

For the $\beta_i^-$, let $\bar{g_i}^{k_i + 1}$ denote reflection about the hyperplane $x_n = 0$. 
The decomposition $\beta = \sum_i \beta_i^+ + \sum_i \beta_i^-$, coupled
with the maps $\bar{g}_i^1 ... \bar{g}_i^{k_i}$ for the $\beta_i^+$, and the maps $\bar{g}_i^1
... \bar{g}_i^{k_i+1}$ for the $\beta_i^-$, gives a decomposition of $\beta$ satisfying hypotheses
1), 2), and 3) of Lemma 2.2. Write $\zeta_i^+ = \bar{g}_i^1 \circ ... \circ \bar{g}_i^{k_i}$ and 
$\zeta_i^- = \bar{g}_i^1 \circ ... \circ \bar{g}_i^{k_i+1}$ respectively, and let $D_i^+ = 
\{x: \zeta_i^+(x) > 0\}$ and $D_i^- = \{x: \zeta_i^-(x) > 0\}$. Let $\bar{D_i}^+$ and $\bar{D_i}^-$
respectively denote their closures. We will see that for each 
$w \in \bar{D_i}^+ \cap (\zeta_i^+)^{-1}(0) = \bar{D_i}^- \cap (\zeta_i^-)^{-1}(0)$, the functions 
$f \circ \zeta_i^+ (x + w)$ and $f \circ \zeta_i^- (x + w)$ each either satisfies the induction
hypothesis or is of the form $(2.13)$. 
Once we establish this, if we assume Theorem 2.4 then hypothesis 4) of Lemma 2.2 holds as well
with $q_w$ the identity map. 
As a result, we can apply Lemma 2.2 and conclude the Main Theorem holds for $f$. Since 
$f$ is completely arbitrary with ${\partial^{\kappa} f \over \partial x_n^{\kappa}} \neq 0$, 
showing $(2.13)$ reduces the Main Theorem to proving Theorem 2.4. 

So we focus our attention on establishing that either $(2.13)$ or the induction hypothesis holds.
We restrict our attention to the 
$f \circ \zeta_i^+ (x + w)$ since the $f \circ \zeta_i^- (x + w)$ are done in an entirely 
analogous fashion. By
definition of the $g_i^j$, equation $(2.14)$ becomes
$$f \circ \zeta_i^+ (x) = \bar{f}_{\kappa}(x_1,...,x_n)x_n^{\kappa} + \sum_{l < \kappa}
\bar{m}_l(x_1,...,x_{n-1})\bar{s}_l(x_1,...,x_{n-1})x_n^l \eqno (2.15)$$
In each (nonzero) term, $\bar{m}_l$ is a monomial, and $\bar{f}$ and the $\bar{s}_l$ are nonvanishing. For each
$w \in supp(\beta_i^+ \circ \zeta_i^+)$ with $w_n = 0$, let $U_w$ be a neighborhood
of $w$ small enough so that $diam(U_w) < {1 \over 2} \min \{|w_j|: w_j \neq 0\}$. Then on the 
neighborhood $U_w - w$ of the origin, we have $x_j + w_j \sim w_j$ if $w_j \neq 0$, and $x_j + w_j
= x_j$ if $w_j = 0$. As a result, on $U_w - w$, each $\bar{m}_l(x_1+w_1,...,x_{n-1}+w_{n-1})$ can
be written as $\hat{m}_l(x_1,...,x_{n-1})\hat{s}_l(x_1,...,x_{n-1})$, where $\hat{m}_l$ is a 
monomial and where $\hat{s}_l(x_1,...,x_{n-1})$ doesn't vanish on $U_w - w$. As a result, we can
let $\tilde{s}_l(x_1,...,x_{n-1}) = \hat{s}_l(x_1,...,x_{n-1})\bar{s}_l(x_1 + w_1,...,x_{n-1} + 
w_{n-1})$ and
write
$$f \circ \zeta_i^+(x + w) = \tilde{f}_{\kappa}(x_1,...,x_n)x_n^{\kappa} + \sum_{l < \kappa}
\hat{m}_l(x_1,...,x_{n-1})\tilde{s}_l(x_1,...,x_{n-1})x_n^l \eqno (2.16)$$
We next change the notation 
in $(2.16)$. If there is an $l$ such that $\hat{m}_l$ is constant, then let $k'$ be the smallest 
such index and define $\hat{f}_{k'}(x_1,...,x_n)$ by
$$\hat{f}_{k'}(x_1,...,x_n)x_n^{k'} = \tilde{f}_{\kappa}(x_1,...,x_n )x_n^{\kappa} + \sum_{l \geq k'}\hat{m}_l(x_1,...,x_{n-1}) 
\tilde{s}_l(x_1,...,x_{n-1})x_n^l \eqno (2.17)$$
Note that ${\hat{f}}_{k'}(0) \neq 0$. As a result, we can write
$$f \circ \zeta_i^+(x + w) =  \hat{f}_{k'}(x_1,...,x_n)x_n^{k'} + \sum_{l < k'}\hat{m}_l(x_1,...,x_{n-1})
\tilde{s}_l(x_1,...,x_{n-1})x_n^l \eqno (2.18)$$

\noindent Here $k' \leq \kappa$, ${\tilde{f}}_{k'}(0) \neq 0$, each $\hat{m}_l$ is a 
nonconstant monomial, and each $\tilde{s}_l$ satisfies $\tilde{s}_l(0) \neq 0$. Notice the right-hand 
sum may be empty. If $k' < \kappa$, then $f \circ \zeta_i^+(x + w)$ satisfies the 
induction hypothesis. If $k' = \kappa$, then $f \circ \zeta_i^+(x + w)$ is of the form $(2.13)$.
So each $f \circ \zeta_i^+(x + w)$ either satisfies $(2.13)$ or the induction hypothesis.
A very similar argument shows the same for $f \circ \zeta_i^- (x + w)$. As a result, once we 
prove Theorem 2.4,
the Main Theorem holds for each $f \circ \zeta_i^+(x + w)$ and $f \circ \zeta_i^-(x + w)$ when 
$w \in supp(\beta_i^+ \circ \zeta_i^+)$ or $supp(\beta_i^+ \circ \zeta_i^+)$ respectively and
$w_n = 0$. So in particular the Main Theorem holds for $w \in supp(\beta_i^+ \circ \zeta_i^+) \cap 
(\zeta_i^+)^{-1}(0)$ or $w \in supp(\beta_i^- \circ \zeta_i^-) \cap (\zeta_i^-)^{-1}(0)$.
As described above $(2.15)$, Lemma 2.2 then applies and the Main Theorem holds for $f$. Thus the
Main Theorem is reduced to proving Theorem 2.4

\noindent Sections 3 and 4 of this paper develop techniques to prove Theorem 2.4.

\noindent {\bf 3. Defining regions via the Newton polyhedron}

Suppose $g(x)$ is some real-analytic function defined on a neighborhood of the origin satisfying 
the hypotheses of Theorem 2.4. Where $a$ denotes a multiindex $(a_1,...,a_n)$, we Taylor expand 
$g(x)$ about the origin:
$$g(x) = \sum_a c_ax^a \eqno (3.1)$$
\noindent {\bf Definition:} Let $S_a = \{(x_1,...x_n) \in \R^n: x_i 
\geq a_i \hbox { for all } i\}$. The {\it Newton polyhedron} $N(g)$ of $g$ is defined to be 
the convex hull of the $S_a$ for which  $c_a \neq 0$. 

Observe that since each $S_a$ is closed unbounded polyhedron, so is their convex hull $N_g$. 
Often the extreme points of $N(g)$ are referred to as the vertices of $N(g)$. We have the 
following well-known fact about Newton polyhedra:

\noindent {\bf Fact:} The vertices of $N(g)$ consist of finitely many points $a$ for 
which $c_a \neq 0$.

\noindent The fact that any separating hyperplane for $N(g)$ contains at least one its extreme 
points can be translated as follows:

\noindent {\bf Lemma 3.1:} Let $(x_1,...,x_n)$ satisfy $1 > x_i > 0$ for all $i$. Then 
for any $w \in N(g)$ there is a vertex $v$ of $N(g)$ for which $x^w \leq x^v$.

\noindent {\bf Proof:} The equation $x^w \leq x^v$ is equivalent to $\log(x) \cdot w \leq
\log(x) \cdot v$, where $\log(x)$ denotes $(\log(x_1),...,\log(x_n))$. Since the components of 
$\log(x)$ are all negative and $N(g) \subset \{y: y_l \geq 0 \hbox{ for all } l\}$, the lemma 
follows from the fact that the hyperplane 
$-\log(x) \cdot y = e$ that intersects $N(g)$ with minimal $e$ must contain $y = v$ for 
some extreme point $v$ of $N(g)$. 

Divide a small cube $(-\eta,\eta)^n$, for a small $\eta$ to be determined by $N(g)$, into $2^n$ 
subcubes via the coordinate planes $\{x: x_i = 0\}$. In the following arguments we will only consider
the subcube $(0,\eta)^n$ as the other $2^n - 1$ are done similarly.  We will  
subdivide $(0,\eta)^n$ into a finite collection of disjoint open sets whose union is 
$(0,\eta)^n$ up to a set of measure zero. The idea behind the subdivision is as follows. Let 
$E$ denote the collection of vertices and faces (of any dimension) of $N(g)$.
Each element of $E$ will correspond to one of the open sets in the subdivision.
For $x$ in the open set corresponding to some $F \in E$, $x^v$ will be large if $v \in F$,
while $x^v$ will be far smaller for $v \notin F$. 

Denote the set of vertices of the Newton polyhedron $N(g)$ by $v(g)$. For each subset $S$ of 
$v(g)$, let $V_S$ be the convex hull of $S$. For each $0 \leq i \leq n$, let $V_{i1}$...$V_{im_i}$
be an enumeration of those $V_S$ of dimension $i$ that are not properly contained in any other
$V_S$ of dimension $i$. We next inductively define some corresponding
sets $W_{ij}$, starting with the $W_{nj}$, then defining the $W_{n-1,j}$, and so
on. The definition of the $W_{ij}$ requires an increasing collection of constants $1 < C_0 < ...
< C_n$ depending on $N(g)$. Specifically, for constants $A_1, A_2 > 1$ depending on $N(g)$,
the $C_i$ can be any collection of constants satisfying
$$C_0 > A_1\,\,\,\,\,\,\,\,\,\,\,\,\,\,\,C_{i+1} > C_i^{A_2} \hbox {    for all }i$$
Rather than trying to define $A_1$ and $A_2$ in advance, we simply stipulate that they are large enough
that the arguments of this section and section 4 work. 

\noindent Where $x^v$ denotes $x_1^{v_1}
...\,x_n^{v_n}$, define $W_{ij}$ to be the interior of the following set:
$$\{(x_1,...,x_n) \in (0,\eta)^n: \hbox { the } v \in v(g) \hbox { with } x^v
\hbox { maximal is in } v(g) \cap V_{ij},$$
$$C_i^{-1} x^v < x^{v'} < C_i x^v \hbox { for all } v, v' \in v(g) \cap V_{ij}, \hbox { and } $$
$$(x_1,...,x_n) \notin W_{i'j'} \hbox { if } i' > i \hbox { or if } i' = i \hbox { and } 
j' < j \eqno (3.2)$$
Note that every $x$ is in at least one closure $\bar W_{ij}$; Let $v \in v(g)$ maximize $x^v$ and
suppose $V_{0j} = \{v\}$; $x$ will be in $\bar W_{0j}$ if it has not already been selected to be in 
one of the previously defined $W_{ij}$.

There are two facts that encapsulate the most important properties of the $W_{ij}$. The first,
Lemma 3.6 below, is that if $\eta$ is  sufficiently small, depending on $N(g)$, then 
the only nonempty $W_{ij}$ are those corresponding to compact faces (including vertices) of $N(g)$.
The second is given by the following lemma.

\noindent {\bf Lemma 3.2.} 
Let $v(g)$ denote the set of vertices of $N(g)$. There are $A_1, A_2 > 1$ such that if $C_0,...,C_n$
are constants with $C_0 > A_1$ and $C_{i+1} > C_i^{A_2}$ for all $i$, then one can define the $W_{ij}$ 
so that

\noindent {\bf a)} Let $i < n$. If the following two statements hold, then  $x \in W_{ij}$.

\noindent {\bf 1)} If $v \in v(g) \cap F_{ij}$ and  $v' \in v(g) \cap
(F_{ij})^c$ we have $ x^{v'} < C_n^{-1} x^v$. \parskip = 0pt

\noindent {\bf 2)} For all $v, w \in v(g) \cap F_{ij}$ we have $C_i^{-1}x^w < x^v < C_ix^w$. \parskip = 12pt

\noindent {\bf b)} There is a $\delta > 0$ depending on $N(g)$, and not on $A_1$ or $A_2$, such that if 
$x \in W_{ij}$, then the following two statements hold.

\noindent {\bf 1)} If $v \in v(g) \cap F_{ij}$ and  $v' \in v(g) \cap
(F_{ij})^c$ we have $ x^{v'} < C_{i+1}^{-\delta} x^v$. \parskip = 0pt

\noindent {\bf 2)} For all $v, w \in v(g) \cap F_{ij}$ we have $C_i^{-1}x^w < x^v < C_ix^w$. \parskip = 12pt

\noindent {\bf Proof:} We start with a). Assume the assumptions of a) hold. The definition $(3.2)$
tells us that $x \in W_{ij}$ unless it
is in $W_{i'j'}$ for some other $(i',j')$ for which $i' \geq i$. Suppose this were the case. Let
$v' \in v(g) \cap V_{i'j'} \cap (V_{ij})^c$; we know $v'$ exists since $V_{i'j'}$ cannot be 
properly contained in $V_{ij}$. Let $v$ be such that $x^v$ is maximal; $v$ must be in 
$v(g) \cap V_{ij} \cap V_{i'j'}$. By assumption $x^{v'} < C_n^{-1}x^v$; on the other hand
since $v$ and $v'$ are in $V_{i'j'}$ we must have $x^{v'} > C_{i'}^{-1}x^v$. Since the $C_i$ are
increasing, we have a contradiction and we are done with a). 

We proceed to part b). Assume that $x \in W_{ij}$. The second condition holds by definition. So
assume $v \in v(g) \cap V_{ij}$ and  $v' \in v(g) \cap 
(V_{ij})^c$ such that $x^{v'} \geq C_{i+1}^{-\delta}x^v$; we will show that if $\delta$ is small
enough we have a contradiction. Since the $w$ with
$x^w$ maximal is in $V_{ij}$, by $(3.2)$ $x^v$ is within a factor of $C_i$ of this $x^w$ and we have 
$$C_i x^v > x^{v'} \geq C_{i+1}^{-\delta} x^v $$
Since $C_{i+1} > C_i^{A_2}$ we therefore have
$$(C_{i+1})^{{1 \over A_2}} x^v > x^{v'} \geq C_{i+1}^{-\delta} x^v $$
Letting $\delta' = \max(\delta,{1 \over A_2})$, this becomes
$$(C_{i+1})^{\delta'} x^v > x^{v'} \geq C_{i+1}^{-\delta'} x^v \eqno (3.3)$$
Let $V_{i+1j'}$ be generated by $V_{ij}$ and $v'$; this is the largest $i+1$ dimensional convex set 
generated by elements of $v(g)$ that contains $V_{ij}$ and $v'$. If $w$ is any element of 
$v(g) \cap V_{i+1j'}$, there are $w_{l} \in v(g) \cap (V_{ij} \cup \{v'\})$ and constants $c_l$
such that
$$w - v = \sum_l c_l (w_{l} - v)$$
This implies that
$${x^w \over x^v} = \prod_l ({x^{w_{l}} \over x^v})^{c_l} $$
By $(3.2)$ and $(3.3)$, since each $w_{l} \in v(g) \cap V_{ij}$ or $w_{l} = v'$, each factor 
${x^{w_{l}} \over x^v}$ is between $C_{i+1}^{-\delta'}$ and $(C_{i+1})^{\delta'}$. 
Consequently there is some constant $d$ depending only on $N(g)$ such that
$$(C_{i+1})^{-d\delta'} < {x^w \over x^v} < (C_{i+1})^{{d \delta'}} $$
So for any $w$, $w'$ in $v(g) \cap V_{i+1j'}$, since ${x^w \over x^{w'}} = {x^w \over x^v}
{x^v \over x^{w'}}$, we have
$$(C_{i+1})^{-2d\delta'} < {x^w \over x^{w'}} < (C_{i+1})^{{2d \delta'}} \eqno (3.4)$$
As long as $A_2$ was chosen to be greater than $2d$, if one sets $\delta < {1 \over 2d}$, 
then $\delta' = \max(\delta,{1 \over A_2})$ is less than ${1 \over 2d}$ and $(3.4)$ implies
that $x$ satisfies the definition $(3.2)$ for $W_{i+1j'}$, unless it has already been chosen to 
even be in a previously defined $W_{i''j''}$. This contradicts that $x \in W_{ij}$; the definition 
$(3.2)$ implies that  $x$ is not in any  $W_{i+1j'}$ or a previously defined $W_{i''j''}$. Thus the 
proof is complete.

The next sequence of results, leading up to Lemma 3.6, shows that if $\eta$ were chosen sufficiently
small, then $W_{ij}$ intersects $(0,\eta)^n$ if and only if
the associated $V_{ij}$ is a vertex or face of $N(g)$. This will allow us to prove Theorem 2.4 under
the assumption that the only nonempty $W_{ij}$ are those that derive from a vertex or face of $N(g)$. 
The proof is
done through several lemmas, each of which eliminates certain possibilities for $W_{ij}$.

\noindent {\bf Lemma 3.3:} For each $j$ there is an $\epsilon_{nj}$ such that any $x \in 
W_{nj}$ satisfies $|x| > \epsilon_{nj}$. Hence if $\eta$ is sufficiently small, $W_{nj}$ is empty.

\noindent {\bf Proof:} Let $v_0$,...,$v_n \in v(g) \cap V_{nj}$ such that the vectors
$v_1 - v_0$,...,$v_n - v_0$ are linearly independent. Then for the $k$th unit coordinate vector 
${\bf e}_k$, 
we may write ${\bf e}_k = \sum_{l=1}^n c_{jkl}(v_l - v_0)$ for some constants $c_{jkl}$. As a 
result, for each $x \in W_{nj}$, we have
$$x_k = \prod_{l=1}^n ({x^{v_l} \over x^{v_0}})^{c_{jkl}} \eqno (3.5)$$
Definition $(3.2)$ stipulates that each ${x^{v_l} \over x^{v_0}}$ is bounded above and below by a constant.
Therefore $(3.5)$ gives that each $x_k$ is also bounded below by a constant, and we are done.

\noindent {\bf Lemma 3.4:} If $V_{ij}$ intersects the interior of $N(g)$ or the interior of an
unbounded face of $N(g)$, then there is a constant $\delta_{ij}$ such that any $x \in W_{ij}$ 
satisfies $|x| > \delta_{ij}$. Hence in this case too, if $\eta$ were chosen small enough $W_{ij}$
is empty.

\noindent {\bf Proof:} Let $p$ be a point of $V_{ij}$ intersecting the interior of 
$N(g)$ or the interior of an unbounded face of $N(g)$. In either case, denote this interior by $I$.
Because $I$ is unbounded, there is a vector $s$ whose components are all nonnegative with at least
one positive, such that $p - s$ is still in $N(g)$. Because $p - s$ is in $N(g)$, 
$p - s$ is of the form $q + s'$, where $q$ is a convex combination $\sum_l t_l v_l$
of elements of $v(g)$ and where
each component of $s'$ is nonnegative. Letting $r = s + s'$, we have that $p = q + r$, where
each component of $r$ is nonnegative with some component $r_k$ being strictly positive. 

The point $p$ is in $V_{ij}$, so we may write $p$ as a convex combination $\sum_l r_l v_l$
of elements of $V_{ij} \cap v(g)$.  For any $x \in W_{ij}$ we have
$x^p = \prod_l(x^{v_l})^{r_l}$; since $\sum_l r_l = 1$  we have
$$x^p \geq \min_l x^{v_l} > C_i^{-1} \max_{v \in v(g)} x^v \eqno (3.6)$$
The latter inequality follows from the definition $(3.2)$. 
Because the entries of $p - q$ are nonnegative with $p_k - q_k > 0$, we have
$$ x_k^{p_k - q_k}x^q \geq x^p \eqno (3.7)$$
On the other hand we also have
$$x^q = \prod_l(x^{v_l})^{t_l} \leq \max_l x^{v_l} \leq \max_{v \in v(g)} x^v \eqno (3.8)$$
Combining $(3.6) - (3.8)$ we get 
$$C_i^{-1} \max_{v \in v(g)} x^v < x_k^{p_k - q_k} \max_{v \in v(g)} x^v$$
This implies that $x_k \geq C_i^{- {1 \over p_k - q_k}}$ and we are done.

\noindent {\bf Lemma 3.5:} If for some nonempty $W_{ij}$ the set $V_{ij}$ intersects the 
interior of a bounded face $F$ of $N(g)$, then $V_{ij}$ contains $F$.

\noindent {\bf Proof:} The proof is by contradiction. Suppose for a nonempty $W_{ij}$, the
set $V_{ij}$ intersects $F$ at a 
point $p$ in the interior of $F$, but $V_{ij}$ doesn't contain $F$. Since $V_{ij}$ and $F$
are convex but $F $ is not contained in $V_{ij}$, we may let $v' \in v(g) \cap F$ such that $v' \notin 
V_{ij}$. The line starting at $v'$ and passing through $p$ intersects the boundary of $F$ at
a point which we call $q$. There is then $0 < s < 1$ with 
$$ s v' + (1 - s)q = p \eqno (3.9)$$
We rewrite this as
$$q = {1 \over 1 - s}p - {s \over 1 - s}v' \eqno (3.9')$$
We may write $q$ is a convex combination $\sum_l s_l v_l$, where $v_l \in v(g) \cap F$. 
Then for $x \in W_{ij}$, if $v_{max}$ is such that $x^{v_{max}} = \max_{v \in v(g)} x^v$
we have
$$x^q = \prod_l(x^{v_l})^{s_l} \leq x^{v_{max}} \eqno (3.10)$$
Furthermore, by $(3.9')$ we have
$$x^q = (x^p)^{1 \over 1 - s}/ (x^{v'})^{s \over 1 - s} \eqno (3.11)$$
Since $p$ is in $V_{ij}$, $p$ is a convex combination $\sum t_l v_l$ of elements of $V_{ij} 
\cap v(g)$ and by $(3.2)$ we have
$$x^p = \prod_l (x^{v_l})^{t_l} > C_i^{-1}x^{v_{max}} \eqno (3.12)$$
By Lemma 3.2, there is a $\delta > 0$ such that
$$x^{v'} < C_{i+1}^{-\delta} x^{v_{max}} \eqno(3.13)$$
Putting $(3.12)$ and $(3.13)$ into $(3.11)$ we get
$$x^q >  C_i^{-1 \over 1 - s} C_{i+1}^{\delta s \over 1 - s} x^{v_{max}} \eqno (3.14)$$
As long as $A_2$ was chosen larger than ${1 \over \delta s}$, we have 
$$C_{i+1}^{\delta s \over 1 - s} > C_i^{A_2 \delta s \over 1 - s} > C_i^{ 1 \over 1 - s} \eqno (3.15)$$
(Note that the definition of $\delta$ did not depend on $A_2$ so there is no circularity here).
As a result of $(3.15)$ we obtain
$$x^q > x^{v_{max}}$$
This however contradicts $(3.10)$, and we are done.

\noindent We now can prove the following important lemma:

\noindent {\bf Lemma 3.6:} If $\eta$ is sufficiently small, depending on $N(g)$, then 
if $W_{ij}$ is nonempty $V_{ij}$ is a vertex or bounded face of $N(g)$.

\noindent {\bf Proof:} Suppose $W_{ij}$ is nonempty for arbitrarily small
$\eta$. Recall $V_{ij}$ is of dimension $i$. If $i = 0$ there is nothing to prove, so
assume $i > 0$. By Lemma 3.3, $i < n$, and by Lemmas 3.4 and 
3.5, $V_{ij}$ does not intersect the interior of $N(g)$, the interior of an unbounded face of 
$N(g)$, or the interior of any faces of $N(g)$ of dimension greater than $i$. 
But since $V_{ij}$ is $i$-dimensional, we may let $F$ be a (bounded) $i$-dimensional face of $N(g)$
such that $V_{ij}$ intersects the interior of $F$. By Lemma 3.5, $F \subset V_{ij}$. If
$V_{ij}$ contained some point $p$ not on $F$, then since $V_{ij}$ is convex it would 
contain the convex hull of $p$ and $F$, a set of dimension $i + 1$. Since $V_{ij}$ is 
$i$-dimensional, this does not happen. We conclude $F = V_{ij}$ and we are done.

\noindent {\bf Beginning of the proof of Theorem 2.4}

Assume now that $g(x)$ satisfies the hypotheses of Theorem 2.4. In view of Lemma 3.6, in proving 
Theorem 2.4 we may assume that $(0,\eta)^n$ can be written as the union of a set of 
measure zero and the $W_{ij}$ corresponding to vertices and bounded faces of various dimensions of $N(g)$. 
For a given $W_{ij}$, let $F_{ij}$ denote the face or vertex of $N(g)$
for which $V_{ij} = F_{ij}$, and let $e_{ij}$ denote the vertex $\alpha$ of $N(g)$ on $F_{ij}$ whose
$n$th component $\alpha_n$ is maximal; this vertex is unique by the form $(2.13)$. Let $\kappa ' $
denote the $n$th component of $e_{ij}$. So $\kappa ' \leq \kappa$. The following lemma gives upper 
bounds on $g(x)$ and lower bounds on ${\partial^{\kappa ' }g \over \partial x_n^{\kappa ' }}$.
In section 4, each $W_{ij}$ will be subdivided into finitely many $W_{ilp}$, and on each 
$W_{ilp}$ an invertible monomial map will take $W_{ilp}$ to a set $Z_{ilp}$ where the bounds
given by Lemma 3.7 will allow us to use the induction hypothesis on $\kappa$ 
and prove Theorem 2.4.

\noindent {\bf Lemma 3.7:} If the $C_i$ were chosen to increase sufficiently fast, then if 
$\eta$ is sufficiently small there are constants $K, K'$ such for $x \in W_{ij}$ we have
$$ |g(x)| < K x^{e_{ij}}\,\,\,(i > 0) \eqno (3.16a)$$
$$  K' x^{e_{ij}} < |g(x)| < K x^{e_{ij}}\,\,\,(i = 0) \eqno (3.16b)$$
$$ |{\partial^{\kappa ' }g \over \partial x_n^{\kappa ' }}(x)| > K x^{e_{ij}}
x_n^{-\kappa '  } \,\,\,(\hbox {all } i)\eqno (3.16c)$$

\noindent {\bf Proof:} We first prove $(3.16a,b)$. Write $g(x) = \sum_a c_a x^a = \sum_{a \in F_{ij}} 
c_a x^a + \sum_{a \notin F_{ij}} c_a x^a$. In order to prove $(3.16a,b)$ it suffices to show 
two things. First, for some constant $C'$ we will see that
$$\sum_{a \in F_{ij}}|c_a|x^a < C'x^{e_{ij}} \eqno (3.17)$$
Secondly, we will show that given any fixed $\epsilon > 0$, if the $C_i$ are growing fast 
enough and $\eta$ is sufficiently small then we have
$$\sum_{a \notin F_{ij}} |c_a| x^a < \epsilon x^{e_{ij}} \eqno (3.18)$$
Since when $i = 0$ there is only one $a$ in $F_{ij}$, equation $(3.18)$ will automatically
imply the left hand inequality of $(3.16b)$. Equation $(3.16a)$ and the right hand inequality of 
$(3.16b)$ will follow from adding $(3.17)$ and $(3.18)$. 

We consider $(3.17)$ first. If $a \in F_{ij}$, then $a$ can be written as a convex sum 
$\sum_l t_l v_l$ where each $v_l \in v(g) \cap F_{ij}$. consequently we have
$$x^a = \prod_l (x^{v_l})^{t_l} \leq \max_l x^{v_l} < C_ix^{e_{ij}} \eqno (3.19)$$
Adding $(3.19)$ over all $a \in F_{ij}$ gives $(3.17)$. We move to the more difficult $(3.18)$. Every $a$ for
which $c_a x^a$ is nonzero can be written in the form $a = \sum_l t_l v_l + p$, where 
each $v_l \in v(g)$, $\sum_l t_l = 1$, and $p_k \geq 0$ for all
$k$. Let $q^a \in \R^n$ be the vector with integer coordinates such that each component of 
$q^a - \sum_l t_l v_l$ is in $[0,1)$. Since $a$ has integer coordinates, we can write
$a = q^a + r^a$ where every component $r^a_l$ is still greater than or equal to zero. Writing
$N(q^a) = \{ w: w_{ij} \geq q^a_l \hbox{ for all } l\}$, we have that $a \in N(q^a)$. Note that
there are finitely possibilities for $q_a$ since each $q_a$ has integer coordinates and is 
within distance 1 of the convex hull of the elements of $v(g)$. We have
$$\sum_{a \notin F_{ij}} |c_a| x^a \leq \sum_q \sum_{a \in N(q),\,\,a \notin F_{ij}}|c_a| x^a 
\eqno (3.20)$$
In $(3.20)$ we of course only add over the finitely many $q$ that are of the form $q^a$ above.
We divide the sum $(3.20)$ into three parts, depending on where $q$ comes from. Let $L_1$ 
denote the points on $F_{ij}$ with integral coordinates. Let $L_2$ denote the points in 
the convex hull of the elements of $v(g)$ with integral coordinates that are not on $F_{ij}$.
Let $L_3$ denote the remaining possibilities for $q$, namely points not in the convex hull of
the elements of $v(g)$ (but which are within distance 1 of these elements). We have the following.
$$\sum_q \sum_{a \in N(q), a \notin F_{ij}}|c_a| x^a \leq \sum_{q \in L_1} \sum_{a \in N(q)
- {q}}|c_a| x^a + \sum_{q \in L_2}\sum_{a \in N(q)}|c_a| x^a + \sum_{q \in L_3}\sum_{a \in N(q)}
|c_a| x^a \eqno (3.21)$$
We will bound each of the three sums in $(3.21)$; this will give us the desired estimates on
$(3.18)$. First, observe that for $q \in L_1$
$$\sum_{a \in N(q)- {q}}|c_a| x^a < C \sum_l |x_l| x^q < CC_i \sum_l |x_l| x^{e_{ij}}
\eqno (3.22)$$
The last inequality follows from $(3.2)$ since  $q$ and $e_{ij}$ are both on $F_{ij}$ and
$x \in W_{ij}$.
By assuming $\eta$ is sufficiently small, since $\sum_l |x_l| < n\eta$, the right hand side of 
$(3.22)$ can be made less than $\mu x^{e_{ij}}$ for any $\mu$ we'd like.
Moving on to $L_2$, observe that for $q \in L_2$, by continuity of real-analytic functions, if $|x|$ is 
sufficiently small we have 
$$\sum_{a \in N(q)}|c_a| x^a < (|c_q| + 1) x^q \eqno (3.23)$$
Because $q \in L_2$, we can write $q = \sum_l t_l v_l$ for $v_l \in v(g)$, such that at least
one $v_l$ with nonzero $t_l$, say $v_0$, is not on $F_{ij}$. As a result, using Lemma 3.2 and
equation $(3.2)$ we have
$$x^q = \prod_l (x^{v_l})^{t_l} = (x^{v_0})^{t_0}\prod_{l > 0} (x^{v_l})^{t_l} 
\leq (x^{v_0})^{t_0} C_i^{1-t_0}(x^{e_{ij}})^{1 - t_0}$$
$$ < C_i^{-pt_0}(x^{e_{ij}})^{ t_0} \times C_i^{1-t_0}(x^{e_{ij}})^{1 - t_0} = 
C_i^{-pt_0 + 1 - t_0}x^{e_{ij}}\eqno (3.24)$$
Recall $p$ is a positive integer that we may freely choose which determines how fast the 
$C_i$ must grow. For any fixed $\mu$, we can choose $p$ to ensure the right hand side
of $(3.24)$ is at most ${\mu \over |c_q| + 1}x^{e_{ij}}$. This ensures that the right hand side
of $(3.23)$ is at most $\mu x^q$. Next, we move to the terms of $(3.21)$ for $q \in L_3$. For
such $q$, if $\eta$ is sufficiently small then again $(3.23)$ holds. Since $q \in N(g)$ is not a 
convex combination of elements of $v(g)$, we can select a $q'$ which is a convex combination of 
members of $v(g)$ such that each component of $q - q'$ is nonnegative, with at least one 
component, say $q_r - q_r'$, strictly positive. So we have
$$x^q \leq x_r^{q_r - {q_r}'} x^{q'} < C_i x_r^{q_r - {q_r}'} x^{e_{ij}} \eqno (3.25)$$
The last inequality follows from $(3.2)$. If $\eta$ is sufficiently small, we can make 
$C_i x_r^{q_r - {q_r}'} < {\mu \over |c_q| + 1}$ for any $\mu$ one likes, giving
$$\sum_{a \in N(q)}|c_a| x^a < (|c_q| + 1) x^q <  \mu x^{e_{ij}}$$
This gives the desired estimates for a term of $(3.21)$ for $q \in L_3$. So we have now seen 
that each 
$\sum_{a \in N(q)}|c_a|x^a$ or $\sum_{a \in N(q)
- {q}}|c_a| x^a$ in $(3.21)$ can be made less than $\mu x^{e_{ij}}$ for any prechosen $\mu$.
Consequently, the entire sum $(3.21)$ can be made less than any $\epsilon  x^{e_{ij}}$ for
any prechosen $\epsilon$. This gives $(3.18)$ and we are done with part a) of this lemma.

The proof of part c) is quite similar to that of parts a) and b). The Newton polyhedron of 
${\partial^{\kappa ' }g \over 
\partial x_n^{\kappa ' }}$ is obtained by taking the portion of the Newton polyhedron of 
$g$ with ``height'' at least $\kappa ' $ and shifting it downward by $\kappa ' $ units. There is 
exactly one vertex of $F_{ij}$ at height at least $\kappa ' $, namely $e_{ij}$, so the face of 
the Newton polyhedron of ${\partial^{\kappa ' }g \over \partial x_n^{\kappa ' }}$ corresponding
to $F_{ij}$, call it $F_{ij}'$, consists of the single vertex $e_{ij} - (0,...,\kappa ' )$.
Suppose $v'$ is some vertex of the Newton polyhedron of ${\partial^{\kappa ' }g \over 
\partial x_n^{\kappa ' }}$ other than $e_{ij} - (0,...,\kappa ' )$. Then using the form
$(2.3)$, $v' + (0,...,\kappa ' )$ must be in $v(g) \cap (F_{ij})^c$. So for $x \in W_{ij}$ we have
$${x^{e_{ij} - (0,...,\kappa ' )} \over x^{v'}} = {x^{e_{ij}} \over x^{v' + 
(0,...,\kappa ' )}} > C_i^p$$
Consequently, if we write
$${\partial^{\kappa ' }g \over \partial x_n^{\kappa ' }}(x) = d_{e_{ij}}x^{e_{ij}}x_n^{-\kappa ' } + 
\sum_{a \notin F_{ij}'} d_ax^a \eqno (3.26)$$
Then exactly as in the proof of part a), if $\eta$ is small enough and the $C_i$ were chosen to be 
increasing fast enough, the sum
$\sum_{a \notin F_{ij}'} |d_a|x^a$ can be made less than $\epsilon x^{e_{ij}}x_n^{-\kappa ' }$ 
for any $\epsilon$ that we would like. As a result, shrinking $\eta$ if necessary, we can assume
$$ |{\partial^{\kappa ' }g \over \partial x_n^{\kappa ' }}(x)| > {|d_{e_{ij}}| \over 2} x^{e_{ij}}
x_n^{-\kappa ' } \eqno (3.27)$$
This gives part c) of the lemma and we are done.

\noindent {\bf 4: Subdividing the $W_{ij}$, finishing the proof of Theorem 2.4}

In this section we subdivide each $W_{ij}$, modulo a set of measure zero, into finitely many pieces 
$W_{ijp}$. On each  $W_{ijp}$ we will define an invertible monomial map that takes $W_{ijp}$ 
bijectively to a set $Z_{ilp}$. 
An application of Lemma 2.2 on the transformed function will then allow us to use the induction 
hypothesis on $\kappa$. As a result, Theorem 2.4, and 
therefore the Main Theorem, will follow. It is important that after some appropriate reflections each 
$Z_{ilp}$ contains a cube $(0,\rho_{ilp})^n$ and is 
contained in a cube $(0,\rho_{ilp}')^n$ for some $0 < \rho_{ilp} < \rho_{ilp}'$. Hence we introduce
the following definition:

{\bf \noindent Definition:} A set $Q$ is called a {\it positive curved quadrant} if there are
$0 < \rho < \rho'$ such that 
$$(0,\rho)^n \subset Q \subset (0,\rho')^n$$

For each $i$ and $j$ let $f_{ij}$ be the vertex $(f_{ij1},...,f_{ijn})$ on $F_{ij}$ 
such that the component $f_{ijn}$ is minimal; there exists a unique
such vertex by $(2.3)$. Since the face $F_{ij}$ is of dimension $i$, we may let $\{P_l\}_{l=1}^{n-i}$
be separating hyperplanes for $N(g)$ such that $F_{ij} = \cap_{l=1}^{n-i} P_l$. We write these 
hyperplanes as 
$$P_l = \{x: a^l \cdot  x = c^l\}$$
We can assume the $a^l$ have rational coefficients. The hyperplanes satisfy
$$N(g) \subset \cap_{l=1}^{n-i} \{x: a^l \cdot x \geq c^l\} \eqno (4.1)$$
Since $\cap_{m=1}^n\{x: x_m  \geq f_{ijm}\} \subset N(g)$, we also have
$$\cap_{m=1}^n\{x: x_m  \geq f_{ijm}\} \subset \cap_{l=1}^{n-i} \{x: a^l \cdot x \geq c^l\} \eqno 
(4.2)$$
Since $a^l \cdot f_{ij} = c^l$ for all $l$, if we shift $x$ in $(4.1)$ by $-f_{ij}$ we get
$$\cap_{m=1}^n\{x: x_m  \geq 0 \} \subset \cap_{l=1}^{n-i} \{x: a^l \cdot x \geq 0\} \eqno 
(4.3)$$
In the case where $i > 0$, we would like to extend the hyperplanes $a^l \cdot x = 0$ to a collection of 
$n$ independent hyperplanes such that 
$$ \cap_{m=1}^n\{x: x_m  \geq 0\} \subset \cap_{l=1}^n \{x: a^l \cdot x \geq 0\} 
\eqno (4.4)$$
(Note that $(4.4)$ automatically holds when $i = 0$.)
We do the extension for $i > 0$ as follows. The point $(0,...,0,1)$ is not in the span of the 
$a^l$ since by $(2.3)$ each extreme point of $F_{ij}$ must have a different $n$th coordinate. 
So we may define $a^n = (0,....,0,1)$ and 
the vectors $a^1,...,a^{n-i}$ and $a^n$ are linearly independent. We similarly define any remaining
$a^l$ for $i < l < n$ to be unit coordinate vectors such that $a^1,...,a^n$ are linearly
independent. Note that we have
$$ \cap_{m=1}^n\{x: x_m  \geq 0\} \subset \cap_{l=n-i+1}^{n}\{x: a^l \cdot x \geq 0\} 
\eqno (4.5)$$
Combining with $(4.3)$ shows that $(4.4)$ holds.

Since the $a^l \cdot x \geq 0$ are $n$ independent hyperplanes intersecting at the origin, any $n-1$ 
of the hyperplanes 
intersect along a line through the origin. Write the directions of these lines as $b_l$, chosen so
that the $b_l$ have rational components and $a_l \cdot b_l > 0$. The $b_l$ span $\R^n$, so we may write
the $m$th unit coordinate vector ${\bf e}_m$ in the form
$${\bf e}_m = \sum_l d_{lm} b_l \eqno (4.6)$$
\noindent {\bf Lemma 4.1:} The coefficients $d_{lm}$ are all nonnegative rational numbers.

\noindent {\bf Proof:} 
By definition of $b_l$, we have
$$ \cap_{l=1}^n \{x: a^l \cdot x \geq 0\} = \{s: s = \sum_p s_p b_p \hbox { with } s_p \geq 0 \} 
\eqno (4.7)$$ 
Since each ${\bf e}_m$ is in $\cap_{m=1}^n\{x: x_m  \geq 0\} \subset \cap_{l=1}^n \{x: a^l
\cdot x \geq 0\}$, $(4.7)$ says that each $d_{lm}$ 
is nonnegative. Elementary linear algebra gives a formula for the $d_{lm}$ which shows that they are
rational. This completes the proof.

We now do a coordinate change on each $W_{ij}$ for $i > 0$. Denoting the original coordinates of a point $x$ 
by $(x_1,...,x_n)$, we let the new coordinates be denoted by $(y_1,...,y_n)$, where 
$$y_m = \prod_{l=1}^n x_l^{d_{lm}} \eqno (4.8)$$
Although the exponents in $(4.8)$ are not necessarily integers and therefore the coordinate
change is not an invertible monomial map, in this section what we will do is compose two
coordinate changes of the form $(4.8)$ with a map $(z_1,...,z_n) \rightarrow (z_1^N,...,
z_n^N)$ for a sufficiently large $N$; this will ensure the resulting composition is an 
invertible monomial map and thus satisfies the requirements of the Main Theorem. 

Observe that a monomial $x^{\alpha}$ becomes $y^{L(\alpha)}$ in the new coordinates, where $L$
is the linear map such that $L(b_l) = {\bf e}_l$ for all $l$. If $\bar{f}_{ij} = 
(\bar{f}_{ij1},...,\bar{f}_{ijn})$ denotes $L(f_{ij})$, then each $\bar{f}_{ijk} \geq 0$ since
each $d_{lm}$ is nonnegative. Furthermore, $L$ takes each hyperplane $P_l$ to
$\{y: y_l = \bar{f}_{ijl}\}$. Notice that each vertex $v$ of $N(g)$ on $F_{ij}$ is on $P_l$ for $l \leq n - i$. This means that
the $l$th component of $L(v)$ is equal to $\bar{f}_{ijl}$ for $l \leq n - i$. 
So if $v$ and $v'$ are vertices of $N(g)$ on $F_{ij}$, the first $n-i$ components of $L(v - v')$ 
are zero. Hence  ${y^{L(v)} \over y^{L(v')}}$ is a function of
the last $i$ $y$-variables only. Write $y = (s,t)$, where $s$ is the first $n-i$ variables and $t$ 
is the last $i$ variables. Similarly, write $L= (L_1,L_2)$, where $L_1$ is the first $n-i$ 
components and $L_2$ is 
the last $i$ components. Recall from $(3.2)$ that for any such $v$ and $v'$, any  
$x \in W_{ij}$ satisfies the inequalities
$$C_i^{-1} < {x^v \over x^{v'}} < C_i \eqno (4.9a)$$
In terms of the $t$ variables this translates as 
$$C_i^{-1} < {t^{L_2(v)} \over t^{L_2(v')}} < C_i \eqno (4.9b)$$
Write $\log(t) = (\log(t_1),\log(t_2),...,\log(t_n))$. Equation $(4.9b)$ becomes
$$-\log(C_i) < \log(t) \cdot L_2(v - v') < \log(C_i) \eqno (4.10)$$
Since the set of all possible $L_2(v - v')$ for $v$ and $v'$ vertices of $g$ on $F_{ij}$ spans an 
$i$-dimensional space, and since $\log(t)$ is an $i$-dimensional vector, there
must be a constant $d$ depending on the function $g$ such that for each $l$ we have
$$-d\log(C_i) < \log(t_l) < d\log(C_i) \eqno (4.11a)$$
Equation $(4.11a)$ is equivalent to 
$$C_i^{-d} < t_l < C_i^d \eqno (4.11b)$$
In particular, the variables $t_l$ are bounded away from 0.  

Next, observe that $L$ takes the hyperplane $P_n = \{\alpha: \alpha_n = 
f_{ijn}\}$ to the hyperplane $P_n' = \{\alpha: \alpha_n = \bar{f}_{ijn}\}$. Replacing the 
vectors $b_l$ by $cb_l$ for an appropriate positive constant $c$, we may assume that $f_{ijn} = 
\bar{f}_{ijn}$. Since $L$ is linear, this implies that $L$ takes any hyperplane $\{\alpha: \alpha_n = C\}$ to 
itself. This fact is useful for finding expressions analogous to $(3.16a),(3.16c)$ in the $y$ 
coordinates. Let $T$ denote the map of the coordinate change from $y$ to $x$ coordinates, and let
$\tilde{g}(y) = g \circ T(y)$. Then $(3.16a)$ gives 
$$ |\tilde{g}(y)| < K y^{L(e_{ij})} \eqno (4.12a)$$
For the derivatives, we use the chain rule. We have
$${\partial \tilde{g} \over \partial y_n}(y) = \nabla g(Ty) \,\,DT(y)\, {\bf e}_n$$
Here $DT(y)$ denotes the derivative matrix of $T$ at $y$. Note that $DT(y) \,{\bf e}_n$ is the last 
column of $DT(y)$. Since $L$ takes each hyperplane $\{\alpha: \alpha_n = C\}$ to itself, each of 
the functions $x_1$,...,$x_{n-1}$ is a function of the $y_1$,..., $y_{n-1}$ variables only, and
$x_n$ is of the form $\tilde{y}^{\alpha}y_n$ where $\tilde{y} = (y_1,...,y_{n-1})$. Consequently,
for $l < n$ we have
${\partial x_l \over \partial y_n} = 0$, while ${\partial x_n \over \partial y_n} =  
\tilde{y}^{\alpha}$. Hence $DT(y) \,{\bf e}_n = \tilde{y}^{\alpha}{\bf e}_n$, and 
$${\partial \tilde{g} \over \partial y_n} (y) = \tilde{y}^{\alpha}{\partial g \over \partial x_n} 
(Ty)$$
Repeating this $\kappa ' $ times, where $\kappa ' - 1$ is as in Lemma 3.7, we have 
$${\partial^{\kappa ' } \tilde{g} \over \partial y_n^{\kappa ' }} (y) = (\tilde{y}^{\alpha})^{\kappa ' }
{\partial^{\kappa ' } g \over \partial x_n^{\kappa ' }} (Ty)$$
Putting this in $(3.16c)$, we have
$$|{\partial^{\kappa ' } \tilde{g} \over \partial y_n^{\kappa ' }} (y)| > K (\tilde{y}^{\alpha})^{\kappa ' }
x^{e_{ij}}x_n^{-\kappa '  } = K (\tilde{y}^{\alpha})^{\kappa ' }y^{L(e_{ij})}(y_n \tilde{y}^{\alpha})^{-\kappa '  }$$
$$ = K y^{L(e_{ij})}y_n^{-\kappa '  }$$
But the variable $y_n$ is bounded below by $(4.11b)$, so the last equation implies
$$|{\partial^{\kappa ' } \tilde{g} \over \partial y_n^{\kappa ' }} (y)| > K'y^{L(e_{ij})} \eqno (4.12b)$$
This is the inequality we seek. Note that the right hand sides of $(4.12a)$ and $(4.12b)$ are 
the same up to a constant. After doing further coordinate changes in the $s$ variables only
(which do not change $(4.12a)-(4.12b)$), we will be able to factor out a $y^{L(e_{ij})}$ from
$\tilde{g}$, generally resulting in a bounded function with a $\kappa ' $st derivative bounded below.
After an application of Lemma 2.2 we will be able to invoke the induction hypothesis. 
As a result Theorem 2.4, and thus the Main Theorem, will be proved. 

Next, continuing to focus on the $i > 0$ case, we examine how the other inequalities in $W_{ij}$'s 
definition behave under this 
coordinate change. It turns out that the relevant inequalities are those provided by Lemma 3.2.
This lemma says that if $x \in W_{ij}$, $w$ is in the vertex set $v(g)$ of $N(g)$ and on the
face $F_{ij}$, and $w' \in v(g)$ but $w' \notin F_{ij}$, then we have
$$x^{w' - w} < (C_{i+1})^{-\delta}$$
Writing in $y$ coordinates, this becomes
$$y^{L(w' - w)}< (C_{i+1})^{-\delta} \eqno (4.13a)$$
We would like to encapsulate the condition that $x \in (0,\eta)^n$ through an equation analogous to
$(4.13a)$. Shrinking $\eta$ if necessary, we can assume that for each $m$, $x_m = x^{{\bf e}_m} < 
(C_{i+1})^{-\delta}$, and we express this in $y$ coordinates as
$$y^{L({\bf e}_m)}< (C_{i+1})^{-\delta} \eqno (4.13b)$$
Writing $L = (L_1,L_2)$ and $y = (s,t)$ like before,  equations $(4.13)$ become
$$s^{L_1(w' - w)} < (C_{i+1})^{-\delta} t^{L_2(w - w')} \eqno (4.14a)$$
$$s^{L_1({\bf e}_m)} < (C_{i+1})^{-\delta} t^{L_2(-{\bf e}_m)} \eqno (4.14b)$$
Equation $(4.11b)$ says that each component of $t$ 
is between $C_i^{-d}$ and $C_i^d$. So there is a constant $d'$ depending only $N(g)$ such that in
$(4.14)$ one has
$$C_i^{-d'} < t^{L_2(w - w')} < C_i^{d'} \eqno (4.15a)$$
$$C_i^{-d'} < t^{L_2(-{\bf e}_m)} < C_i^{d'} \eqno (4.15b)$$
So as long as $A_2$ from the beginning of section 3 is sufficiently large, equations $(4.14)$ give
$$s^{L_1(w' - w)} < 1 \eqno (4.16a)$$
$$s^{L_1({\bf e}_m)} < 1 \eqno (4.16b)$$
Summarizing, if $x \in W_{ij}$, then the corresponding $(s,t)$ in $y$ coordinates satisfy $(4.9b)$ and
$(4.16a)-(4.16b)$.
We now use in a similar fashion the other inequalities of Lemma 3.2. Namely, $x \in (0,\eta)^n$ is in 
$W_{ij}$ if $(4.9)$ holds and $x$ satisfies the following for all $w \in v(g) \cap F_{ij}$, 
$w' \in v(g) \cap (F_{ij})^c$
$$x^{w'} < C_n^{-1}x^w \eqno (4.17a)$$
Analogous to above, we incorporate the condition $x \in (0,\eta)^n$ by stipulating that $\eta < (C_n)^
{-1}$ and write
$$x^{{\bf e}_m} < C_n^{-1} \eqno (4.17b)$$
Analogous to $(4.14)$, these can be written as 
$$s^{L_1(w' - w)} < (C_{n})^{-1} t^{L_2(w - w')} \eqno (4.18a)$$
$$s^{L_1({\bf e}_m)} < (C_{n})^{-1} t^{L_2(-{\bf e}_m)} \eqno (4.18b)$$
Again using $(4.15)$, there is some $\mu$ such that equations $(4.18)$ hold whenever for all $w' - w$
and all ${\bf e}_m$ we have
$$s^{L_1(w' - w)} < \mu  \eqno (4.19a)$$ 
$$s^{L_1({\bf e}_m)} < \mu  \eqno (4.19b)$$
Hence if a point $(s,t)$ is such that $s$ satisfies $(4.19a)-(4.19b)$ and $t$ satisfies $(4.9a)$, then
the corresponding $x$ is in $W_{ij}$. Putting $(4.16)$ and $(4.19)$ together, 
let $Y_{ij}$ denote the set $W_{ij}$ in the $y$ coordinates. Let $u_1$, $u_2$,... be an enumeration
of the set of all $L_1(w' - w)$ for vertices $w \in F_{ij}$ and $w' \notin F_{ij}$, as well as the 
distinct $L_1({\bf e}_m)$. 
We define the sets $E_1$ and $E_2$ by
$$E_1 = \{ s: 0 < s^{u_l} < \mu  \hbox { for all } l \} \times D_{ij} \eqno (4.20a)$$
$$E_2 = \{ s: 0 < s^{u_l} < 1  \hbox { for all } l \} \times D_{ij} \eqno (4.20b)$$
Then by $(4.16)$ and $(4.19)$ we have 
$$E_1 \subset Y_{ij} \subset E_2 \eqno (4.20c)$$
It is worth pointing out that none of the $u_l$ are zero: If some $\bar{w}_l - \bar{w}_0$ were 
zero this would imply that they came from a $w \in F_{ij}$ and a $w' \notin F_{ij}$ such that
$w' - w$ is a function of only the $t$-variables. This would mean that $w' - w$ is tangent to
$F_{ij}$, which can never happen when $w \in F_{ij}$ and $w' \notin F_{ij}$. If some 
$L_1({\bf e}_m)$ were zero, that would imply ${\bf e}_m$
is a function of the $t$ variables only, meaning that ${\bf e}_m$ is tangent to $F_{ij}$. Since 
$F_{ij}$ is a bounded face, this cannot happen either.

Equations $(4.20a)-(4.20c)$ are for $i > 0$, and there are analogous equations when $i = 0$.
Fortunately, these require less effort to deduce; a coordinate change is not required.
There is a single vertex $v$ on a given $F_{0j}$.
Lemma 3.2 tells us that if $\mu$ is sufficiently small, if we define 
$$F_1 = \{ x \in (0,\eta)^n: x^{v'} < \mu x^v \hbox { for all } v' \in v(g) - \{v\}\}$$
$$F_2 = \{ x \in (0,\eta)^n: x^{v'}< x^v \hbox { for all } v' \in v(g) - \{v\}\}$$
Then we have $F_1 \subset W_{0j} \subset F_2$. To combine this with the $i > 0$ case, we
rename the $x$ variables $s$ and define $Y_{0j} = 
W_{0j}$. Let $\{u_l\}_{l > 0} $ be an enumeration of the $v' - v$ for
$v' \in v(g) - \{v\}$ as well as the unit coordinate vectors ${\bf e}_m$. When $i = 0$ define 
$$E_1 = \{ s:0 < s^{u_l} < \mu  \hbox { for all } l > 0\}$$
$$E_2 = \{ s:0 < s^{u_l} <  1 \hbox { for all } l > 0\}\eqno(4.21)$$
Then, shrinking $\mu$ to less than $\eta$ if necessary, like above we have $E_1 \subset Y_{0j} \subset 
E_2$.

In the remainder of this section, we consider the $i > 0$ and $i = 0$ cases together. We still have some
work to do. Namely, we would like to replace the sets $\{ s: 0 < s^{u_l} < \mu  \hbox { for all } l \}$
or $\{ s: 0 < s^{u_l} < 1 \hbox { for all } l \}$ by cubes. To this end, 
we will divide up $Y_{ij}$ in the $s$ variables into finitely many pieces. A 
coordinate change in the $s$ variables will be performed on each piece taking it to a set which is a 
positive curved quadrant. This is done as follows. For $i > 0$ let $E_1'$ and $E_2'$ be defined by
$$ E_1' = \{ s :0 < s^{u_l} < \mu \hbox { for all } l > 0\}$$
$$ E_2' = \{ s :0 < s^{u_l} < 1 \hbox { for all } l > 0\}$$
When $i = 0$, let $E_1' = E_1$ and $E_2' = E_2$.
Writing  $S = (S_1,..,S_{n-i}) = (\log(s_1),..,\log(s_{n-i}))$, in the $S$ coordinates $E_2'$ becomes
the set $E_2^S$ given by
$$E_2^S = \{ S: S \cdot u_l < 0  \hbox { for all } l\} $$
The set of $S$ satisfying $(4.21)$ is the intersection of several hyperplanes passing through 
the origin. We subdivide $E_2^S$ via the $n-i$ hyperplanes $S_m = 0$, resulting in (at most)
$2^{n-i}$ pieces which we call $E_2^{S,1}$, $E_2^{S,2}$,... We focus our attention on the one for which
all $S_m > 0$, which we assume is $E_2^{S,1}$. The intersection of $E_2^{S,1}$ with the hyperplane
$\sum_m S_m = 1$ is a polyhedron, which we can triangulate into finitely simplices $\{Q_p\}$ whose 
vertices all have rational coordinates. By taking the convex hull of these $Q_p$'s with the origin,
one obtains a triangulation of $E_2^{S,1}$ into unbounded $n$-dimensional
``simplices'' which we denote by $\{R_p\}$. Each $R_p$ has $n$ unbounded faces of dimension $n-1$
 containing the origin. The equation for a given face can be written as $S \cdot q^{p,l} = 0$, where
each $q^{p,l}$ has rational coordinates, so that
$$R_p = \{S : S \cdot q^{p,l} < 0 \hbox { for all } 1 \leq l \leq n - i \} \eqno (4.22)$$
Hence $\cup R_p = E_2^{S,1}$. The other $E_2^{S,m}$ can be similarly subdivided. We 
combine all simplices from all the $E_2^{S,m}$ into one list $\{R_p\}$. 
Note each $R_p$ on the combined list satisfies $(4.22)$. Furthermore, the $R_p$ are disjoint and
 up to a set of measure zero $E_2^S = \cup_p R_p$. Converting back now into $s$
coordinates, for $i > 0$ we define 
$$Y_{ijp} = \{(s,t) \in Y_{ij}: \log(s) \in R_p\} = \{(s,t) \in Y_{ij}:0 < s^{q^{p,l}} < 1 \hbox 
{ for all } 1 \leq l \leq n - i \}\eqno (4.23a)$$
When $i = 0$ we let
$$Y_{0jp} = \{s \in Y_{0j}: \log(s) \in R_p\} = \{s \in Y_{0j}:0 < s^{q^{p,l}} < 1 \hbox 
{ for all } 1 \leq l \leq n \}\eqno (4.23b)$$
Then the $Y_{ijp}$ are disjoint and up to a set of measure zero we have
$$\cup_p Y_{ijp} = Y_{ij} \subset  E_2\eqno (4.24)$$
On each $Y_{ijp}$ we shift from $y = (s,t)$  coordinates 
(or $y = s$ coordinates if $i = 0$) to $z = (\sigma,t)$ coordinates (or $z = \sigma$ coordinates if
$i = 0$), where $\sigma$ is  defined by
$$\sigma_l = s^{q^{p,l}} \hbox { for } l \leq n - i\eqno (4.25)$$
In the new coordinates, $Y_{ijp}$ becomes a set $Z_{ijp}$ where 
$$ Z_{ijp} \subset (0,1)^{n-i} \times D_{ij}\,\,\,\,\,\,\,\,\,\,(i > 0) \eqno (4.26a)$$
$$ Z_{ijp} \subset (0,1)^n\,\,\,\,\,\,\,\,\,\,(i = 0) \eqno (4.26b)$$
Let $W_{ijp}$ denote the set $Z_{ijp}$ in the original 
$x$ coordinates. So the $W_{ijp}$ are disjoint open sets and up to a set of measure zero 
$\cup_p W_{ijp} = W_{ij}$.

\noindent {\bf Lemma 4.2.} If $i > 0$, write $z = (\sigma, t)$, where $\sigma$ denotes the first $n - i$
 components and $t$ the last $i$ components. For any vector $w$, we denote by $(w',w'')$ the vector 
such that the
monomial $x^w$ transforms to $\sigma^{w'}t^{w''}$ in the $z$ coordinates. In the case where $i = 0$,
write $z = \sigma$ and say that $x^w$ transforms into $\sigma^{w'}$.

\noindent {\bf a}) If $w$ is either a unit coordinate vector ${\bf e}_l$, or of the form $v' - v$ for 
$v$ a vertex of $g$ in $F_{ij}$ and $v'$ a vertex of $g$ not in $F_{ij}$, then each component of 
$w'$ is nonnegative, with at least one component positive.

\noindent {\bf b}) If each component of $w$ is nonnegative, then so is each component of $w'$ and $w''$.

\noindent {\bf c}) There exists some $\mu' > 0$ such that for all $i,j$, and $p$
$$(0,\mu')^{n-i} \times D_{ij} \subset Z_{ijp} \subset (0,1)^{n-i} \times D_{ij} \,\,\,\,\,\,\,\,\,\,(i > 0)\eqno (4.27a)$$
$$(0,\mu')^n \subset Z_{ijp} \subset (0,1)^n\,\,\,\,\,\,\,\,\,\,(i = 0)\eqno (4.27b)$$
In particular, when $i > 0$, for fixed $t$ the cross-section of $Z_{ijp}$ is 
a positive curved quadrant.

\noindent {\bf Proof.} We assume that $i > 0$; the $i = 0$ case is done exactly the same way. If 
$w$ is of one the forms of part a), then the monomial $x^w$ in the $x$ 
coordinates becomes a monomial of the form $s^{u_m}t^a$ in the $y$ coordinates, where the $u_m$ 
are as before.
Since $Y_{ijp} \subset E_2$, where $E_2$ is as in $(4.20)$ or $(4.21)$, whenever each $s^{q^{p,l}} < 1$ for each $l$
we have $s^{u_m} < 1$ for each $m$. Thus if we write $s^{u_m} = \prod_l (s^{q^{p,l}})^{\alpha_l}$, each
$\alpha_l$ must be nonnegative; otherwise we could fix any $s^{q^{p,l}}$ for which $\alpha_l$ is 
nonnegative, and let the remaining $s^{q^{p,l}}$ go to zero, eventually forcing $s^{u_m} = \prod_l 
(s^{q^{p,l}})^{\alpha_l}$ to be greater than 1. This means that the $\alpha_l$ are nonnegative. If
they were all zero, this would mean $u_m = 0$ which cannot happen by the discussion after $(4.20c)$.
So at least one $\alpha_l$ is positive.
Since $s^{u_m}t^{v}$ transforms into $\sigma^{\alpha_l}t^v$ in the $z$ coordinates, we have part 
a) of this lemma.

Next, we saw that any $x_l$ transforms into some $s^{a_l}t^{b_l}$ in the $y$ coordinates, where each
component of $a_l$ and $b_l$ is nonnegative . When transforming from $x$ to $z$ coordinates, by part
a) $x_l$ transforms into some  $\sigma^{a_l'}t^{b_l}$ with $a_l'$ having nonnegative components. 
Hence part b) holds for the $x_l$. Therefore it holds for any $x^w$ with each
component of $w$ nonnegative.

Moving to part c), the right-hand sides follow from $(4.26)$. As for the left hand sides, from the 
expression $s^{u_m} = 
\prod_l (s^{q^{p,l}})^{\alpha_l}$ with nonnegative $\alpha_l$, there is a $\mu' > 0$ such that each 
$s^{u_m} < \mu$ whenever $s^{q^{p,l}} < \mu'$ for all $l$. So if $s^{q^{p,l}} < \mu'$ for 
each $l$ and $t \in D_{ij}$, then $(s,t) \in E_1$. By $(4.20c)$, we conclude that whenever 
$s^{q^{p,l}} < \mu'$ for all $l$ and if $t \in D_{ij}$, then $y = (s,t)$ is in $Y_{ijp}$. In the
$z$ coordinates this becomes the left hand inequality of $(4.27a)$ for $i > 0$. When $i = 0$, the 
same argument holds; 
whenever $s^{q^{p,l}} < \mu'$ for each $l$ then $s \in E_1$ and $(4.27b)$ follows. Thus we are done 
with the proof of Lemma 4.2.

Lemma 4.2 tells us that $x_m = z^{L_3({\bf e}_m)}$ with each component of each $L_3({\bf e}_m)$ being
nonnegative, but for the $z$ to $x$ coordinate change to satisfy the conditions of the Main 
Theorem we need the components to be integers. This is easy to accomplish. We would like
to replace each $z_l$ by $z_l^{N_l}$ for some large integers $N_l$, and we can do this by 
replacing the definition $z_l = s^{q^{p,l}}$ by $z_l = s^{q^{p,l} / N_l}$ in the above arguments.
Lemma 4.3 still holds (possibly with a different $\mu'$), and the components of each 
$L_3({\bf e}_m)$ are now nonnegative integers. 

Furthermore, we can ensure that $(4.12a)-(4.12b)$ 
still hold by stipulating that $N_n = 1$; the $x$ to $y$ coordinate change takes $x_n$ to some
$\tilde{y}^{\alpha}y_n$ where $\tilde{y} = (y_1,...,y_{n-1})$, and the $z$ to $y$ coordinate 
change is in the first $n-i$ components only. When $i > 0$, we let $T$
be the coordinate change from $z$ to $x$ coordinates and define $\bar{g}(z) = g \circ T(z)$.
Then $(4.12)$ gives the following, where $\kappa ' \leq \kappa $ is as before.
$$ |\bar{g}(z)| < K z^{L_3(e_{ij})} \eqno (4.28a)$$
$$|{\partial^{\kappa ' } \bar{g} \over \partial z_n^{\kappa ' }} (z)| > K'z^{L_3(e_{ij})} 
\eqno (4.28b)$$
We split $z = (\bar{z},z')$, where $z'$ are the $t$ variables.
We correspondingly write $L_3 = (L_4,L_5)$. Since by $(4.10)$ the $z'$ variables are bounded above
and below, for some constant $K''$ equations $(4.28a)-(4.28b)$ give
$$ |\bar{g}(z)| < K''{\bar{z}}^{L_4(e_{ij})} \eqno (4.29a)$$
$$|{\partial^{\kappa ' } \bar{g} \over \partial z_n^{\kappa ' }} (z)| > K''{\bar{z}}^{L_4(e_{ij})} 
\eqno (4.29b)$$
Since $\bar{g}$ is defined on a neighborhood of the closure $\bar Z_{ijp}$, $(4.28a)$ implies 
for some real-analytic $h(z)$ the function $\bar{g}(z)$ can be written as
$$\bar{g}(z) = {\bar{z}}^{L_4(e_{ij})}h(z) \eqno (4.30)$$
From $(4.29a)-(4.29b)$, $h(z)$ satisfies
$$|h(z)| < K''' \eqno (4.31a)$$
$$|{\partial^{\kappa ' } h \over \partial z_n^{\kappa ' }} (z)| > K''' \eqno (4.31b)$$
When $i = 0$, one has something even stronger. Equation $(3.16b)$ translates into
$$ K'z^{L_3(e_{ij})} < |\bar{g}(z)| < Kz^{L_3(e_{ij})} \eqno (4.31c)$$
So we may write $\bar{g}(z) = z^{L_3(e_{ij})}h(z)$, where $h(z)$ satisfies
$$ K' < h(z) < K \eqno (4.31d)$$
As a result, when $i = 0$ no more resolving of singularities is needed; $\bar{g}(z)$ is already
a monomial times a nonvanishing function on a positive curved quadrant. Nonetheless, we will 
include the $W_{0jp}$ in the remainder of our arguments so as to have a single unified algorithm.

We are now in a position to complete the proof of Theorem 2.4. We have divided $(0,\eta)^n$ into
the sets $W_{ijp}$ each of which, which after an appropriate invertible monomial map, becomes a
set $Z_{ijp}$ on which $(4.30), (4.31a,b)$ or $(4.30), (4.31c,d)$ hold. To simplify the 
notation, we
let $\{W_r\}$ denote the list of all $W_{ijp}$. Thus each $W_r$ has an invertible monomial map
$\zeta_r$ that 
takes a set $Z_r$ to $W_r$, where $Z_r$ is one of the $Z_{ijp}$. In
particular, $(4.30), (4.31a,b)$ or $(4.30), (4.31c,d)$ holds on each $Z_r$.  

The goal now is
to use these equations along with the induction hypothesis and apply Lemma 2.2 so
as to prove Theorem 2.4 and therefore the Main Theorem. We may write the following disjoint union, 
up to set of measure zero:
$$W_r = \cup_s \cap_t\{x \in E: {p_{rst}(x) \over q_{rst}(x)} < c_{rst}\}\eqno (4.32)$$
Here each $p_{rst}(x)$ and $q_{rst}(x)$ are monomials and each $c_{rst}$ is a positive constant. 

For a small $c > 0$ let $\alpha_c(x)$ be a nonnegative function in $C^{\infty}
(0,\infty)$ such that $\alpha_c(x) = 1$ for $x <  1$ and $\alpha_c(x) = 0$ for $x > 1 + c$.
Observe that one can write $1 - \alpha_c(x)$ as ${\bar{\alpha}}_c({1 \over x})$, where 
${\bar{\alpha}}_c(x)$ is a nonnegative bump function on $(0,\infty)$ equal to 1 for $x < 
{1 \over 1 + c}$ and zero for $x > 1$. Note that 
$$\alpha_c({p_{rst}(x) \over c_{rst}q_{rst}(x)}) + \bar{\alpha_c}
({ c_{rst}q_{rst}(x) \over p_{rst}(x) }) = 1 \eqno (4.33)$$
Consequently, we have
$$\chi_E(x) \prod_{rst} [\alpha_c({p_{rst}(x) \over c_{rst}q_{rst}(x)}) + \bar{\alpha_c}
({ c_{rst}q_{rst}(x) \over p_{rst}(x) })] = \chi_E(x) \eqno (4.34)$$
If $l$ denotes the number of triples $(r,s,t)$ that appear in $(4.34)$, the product $(4.34)$ can be 
written as the sum of $2^l$ terms, each a product of $l$ $\alpha_{c}$ or ${\bar{\alpha}}_{c}$
factors. For a fixed $(r,s)$ the sum of all the terms of this sum which contain an $\alpha_{c}
({p_{rst}(x) \over c_{rst}q_{rst}(x)})$ factor for every $t$ is given by
$$ \prod_t\alpha_c({p_{rst}(x) \over c_{rst}q_{rst}(x)}) \eqno (4.35)$$
Note that $(4.35)$ is equal to 1 on the set $W_{rs}$ defined by
$$W_{rs} = \cap_t\{x \in E: {p_{rst}(x) \over q_{rst}(x)} < c_{rst}\}\eqno (4.36)$$
Hence the sum of $(4.34)$ over $(r,s,t)$ for which $\alpha_{rst}$ appears for every $t$ 
corresponding to at least one $(r,s)$ will be equal to 1 on the union of all $W_{rs}$, which is exactly
$(0,\eta)^n$. Denote this sum by $S^c(x)$. We write out the sum as $S^c(x) = \sum_i 
\beta_i^c(x)$, and each $\beta_i^c(x)$ can be written as a quasibump function in the form
$$\beta_i^c(x) = \prod_{j=1}^{l}\gamma_j^c({r_{ij}(x) \over d_{ij} s_{ij}(x)}) \eqno (4.37)$$
Here each $\gamma_j^c$ is either an $\alpha_c$ or a $\bar{\alpha_c}$, each $r_{ij}$ or $s_{ij}$ is
either a $p_{rst}$ or a $q_{rst}$, and each $d_{ij}$ is either a $c_{ij}$ or a ${1 \over c_{ij}}$.
Note that by $(4.35)$, $\lim_{c \rightarrow 0} \chi_E(x)\,\beta_i^c(x)$
is the characteristic function of a subset of some $W_{rs}$, which in turn is a subset of some $W_r$.
Similarly, if $\zeta_i$ denotes
monomial map taking $Z_r$ coordinates to $W_r$ coordinates, then $\lim_{c \rightarrow 0} 
\chi_E(x)\,\beta_i^c \circ \zeta_i (x)$ becomes the characteristic function of a subset of 
$Z_r$. In particular, we may assume $c$ to be small enough such that equations $(4.30),$ 
$(4.31a,b)$ or $(4.30)$, $(4.31c,d)$ hold on $supp(\chi_E(x)\,\beta_i^c \circ \zeta_i (x))$.

The above considerations were for the cube $(0,\eta)^n$, but clearly the analogous decompositions
can be done on the $2^{n} - 1$ reflected cubes, with the corresponding $\zeta_i$ being some 
reflections followed by an invertible monomial map. We are now in a position to verify that the 
conditions of Lemma 2.2 hold for the function $g(x)$. One puts the $\chi_E(x)\,\beta_i^c(x)$ and the 
corresponding terms for the other $2^n - 1$ cubes into one list $\{\beta_i(x)\}$ and defines 
$\beta(x) = \sum_i\beta_i(x)$. Then this decomposition, using the maps $\zeta_i$, satisfy the 
hypotheses of Lemma 2.2. To see this, we argue as follows:
Everything up to 1) follows directly from the definitions, if we let $D_i 
= \{z: \beta_i \circ \zeta_i(z) > 0\}$. Part 2) holds since $\zeta_i$ is the composition 
of reflections and an invertible monomial map. Part 3) is immediate from $(4.37)$.

As for 4), we break into two cases. Case 1 is when $\kappa '$ of $(4.28)-(4.31)$ satisfies
$\kappa ' \leq \kappa - 1$ for all $i$, and case 2 is when for at least one $i$ we have
$\kappa ' = \kappa$. We examine case 1 first. Since 
$\kappa ' \leq \kappa - 1$, by $(4.31a)$ or $(4.31d)$ we may apply the induction hypothesis on 
$h(z)$ for a given $w$ in 
$supp(\beta_i \circ \zeta_i)$, giving a neighborhood $U_w$ of the origin such that the Main Theorem
holds for 
$h \circ \zeta_i  (z + w)$ on $U_w$. In order to verify that 4) holds and thereby
be able to use Lemma 2.2, we will show that the Main Theorem also holds for $\bar{g} \circ \zeta_i 
 (z + w)$ on a neighborhood of $z = 0$. (We use $q_w = 0$ here). This will follow if we can show that
each $\bar{g} \circ \zeta_i  (\Psi_w(x) + w)$ is the product of a monomial and a nonvanishing 
function, where 
$\Psi_w$ denotes the function from the Main Theorem such that $h \circ \zeta_i (\Psi_w(x) + w)$
is a monomial times a nonvanishing function. Note that by $(4.30)$ or $(4.31)$,
$\bar{g}(z) = m(z)h(z)$ for some monomial $m(z)$. Hence if we can show that for each $w$,
$m \circ \zeta_i  (\Psi_w(x) + w)$ is a monomial times a nonvanishing function, then the
same will hold for $\bar{g}(x)$. This will give part 4) of Lemma 2.2, which in turn 
will complete our proof of Theorem 2.4 for case 1 (By the discussion in section 2, this 
in turn implies the Main Theorem for case 1.)

To achieve this, we first observe that $z + w$ is comparable to a monomial near $z = 0$. This
is true because for any $p$ such that $w_p = 0$, the $p$th component $(z + w)_p$ of $z + p$
is equal to $z_p$, while for any $p$ such that $w_p \neq 0$, $(z + w)_p \sim w_p$ on a 
small enough neighborhood of $z = 0$. So on a neighborhood of $z = 0$, we may write $z + w =
c_0(z)m_0(z)$, where $m_0(z)$ is a monomial and $c(0) \neq 0$. Next,
by 4) of the Main Theorem on $h \circ \zeta_i  (z + w)$, each $p$th component of
$\Psi_w(x)$ can be analogously written as $c_p(x)m_p(x)$ near $x = 0$. 
Since $m \circ \zeta_i$ is also a monomial, in a neighborhood of $x = 0$ we can compose these maps 
and write
$$m \circ \zeta_i  (\Psi_w(x) + w) = c''(x)m''(x) \eqno (4.38)$$
Like before, $m''(x)$ is a monomial and $c''(0) \neq 0$. By the above discussion, this implies that
the last condition 4) of Lemma 2.2 holds for $\bar{g} \circ \zeta_i  (\Psi_w(x) + w)$. 
Thus the proof of Theorem 2.4 is complete in case 1, and therefore the proof of the Main Theorem is 
also complete for case 1.

We now move on to case 2. One can apply the induction hypothesis like in case 1 except for those
$w$ in the support of a $\beta_i \circ \zeta_i$ for which $\kappa ' = \kappa$ and 
${\partial^{\kappa ' - 1} h \over \partial z_n^{\kappa ' - 1}} (w) = 
{\partial^{\kappa - 1} h \over \partial z_n^{\kappa - 1}} (w) = 0$. For such a $w$, because 
${\partial^{\kappa} h \over \partial z_n^{\kappa}} (w) \neq 0$, we may apply the implicit 
function theorem and obtain a surface $S_w$ containing $w$ such that in a 
neighborhood of $w$ the zeroes of the function 
${\partial^{\kappa - 1} h \over \partial z_n^{\kappa - 1}} (z)$ are the points of $S_w$. 
Furthermore, after some quasi-translation $q_w$ in the $z_n$ variable fixing $w$, $S_w$ becomes the 
hyperplane $\{z_n = w_n\}$. 

Consider the function $H_w(z) =  h \circ q_w (z  + w)$. Then $H_w(z)$ falls under case 1; that is, 
since ${\partial^{\kappa - 1} H_w \over \partial z_n^{\kappa - 1}} (z) = 0$ only if $z_n = 0$, the 
decomposition $(2.13)$ on the function $H_w(z)$ in place of $g(x)$ will not have an $l = \kappa - 1$
term. As a result, we can the Main Theorem for case 1 and say that there is a neighborhood $U_w$
of the origin such that $H_w$ satisfies the Main Theorem on $U_w$. Hence condition 4) of Lemma 2.2
holds for the function $h$ at $w$. The proof of case 2 is now completed the way we completed case 1; 
we apply Lemma 2.2 to the function $\bar{g}(z) = m(z)h(z)$, 
this time letting $q_w$ be the identity map except in the special cases above, and 
letting $q_w$ be as above otherwise. (We do not have to worry about $w$ with $w_n = 0$ in the statement
of Lemma 2.2; by $(4.10)$ $|w_n|$ is bounded below). This completes the proof of Theorem 2.4 as well as the Main 
Theorem. 

\noindent {\bf 5. The Proofs of Theorem 1.1 and Lojasiewicz's Inequality}

\noindent We start with the proof of Lojasiewicz's inequality which is the more
straightforward of the two. Let $f_1(x)$ and $f_2(x)$ be real-analytic 
functions defined on a neighborhood $V$ of a compact set $K$ such that
$\{x \in V: f_2(x) = 0\} \subset \{x \in V: f_1(x) = 0\}$. For each $x \in K$ with $f_1(x)f_2(x)
= 0$, let $\phi^x$ 
be a bump function supported in $V$ for which the corollary to the Main Theorem holds for
$f_1$ and $f_2$. If $f_1(x)f_2(x) \neq 0$, let $\phi^x$ be a bump function with $\phi^x(x) \neq 0$
such that $f_1f_2 \neq 0$ on $supp(\phi^x)$. 
By compactness, we can let $\{\phi^{x_j}\}$ be a finite collection of these
functions such that $K \subset \cup_j \{x: \phi^{x_j}(x) > 0\}$. Let $V' = \cup_j 
\{x: \phi^{x_j}(x) > 0\}$; this will be a set on which Lojasiewicz's inequality holds. 

Clearly, it suffices to show Lojasiewicz's inequality on each $supp(\phi^{x_j})$ for which
$f_1(x_j)f_2(x_j) = 0$. Write 
$\phi^{x_j} = \sum_i \phi^{x_j}_i$ as in the Main Theorem. It suffices to show Lojasiewicz's 
inequality on each $supp(\phi^{x_j}_i)$. Let $\Psi_{ij}$ denote the composition of coordinate
changes as in the Main Theorem. It suffices to show Lojasiewicz's inequality for 
$f_1 \circ \Psi_{ij}$ and $f_2 \circ \Psi_{ij}$ on $supp(\phi^{x_j}_i \circ \Psi_{ij})$. By the 
Main Theorem, on $supp(\phi^{x_j}_i \circ \Psi_{ij})$ we may write 
$$f_1 \circ \Psi_{ij}(x) = c_1(x)m_1(x), \hbox{   } f_2 \circ \Psi_{ij}(x) = c_2(x)m_2(x) \eqno 
(5.1)$$
Here $m_1(x)$ is some monomial $\prod_{l=1}^n x_l^{\alpha_l}$, $m_2(x)$ is some monomial
$\prod_{l=1}^n x_l^{\beta_l}$, and $c_1(x), c_2(x)$ are functions that don't vanish on $supp(\phi^{x_j}_i
\circ \Psi_{ij})$. 
In order to show Lojasiewicz's inequality, it suffices to show that if $\beta_l > 0$ then 
$\alpha_l > 0$. We do this by contradiction; suppose $\beta_l > 0$ but 
$\alpha_l = 0$. By the corollary to the Main Theorem, $\Psi_{ij}$ extends to some ball $B_{ij}$ 
centered at the origin.
If $B_{ij}$ is small enough, we have that $\Psi_{ij}(B_{ij}) + x_j \subset V$ and $f_1 \circ 
\Psi_{ij}$ and $f_2 \circ \Psi_{ij}$ satisfy $(5.1)$ on $B_{ij}$ with $c_1(x)$ and $c_2(x)$ 
nonvanishing. Let
$z \in B_{ij}$ such that $z_l = 0$ but $z_{l'} \neq 0$ for $l' \neq l$. Then since $\beta_l > 0$, 
we must have $f_2 \circ \Psi_{ij}(z) = 0$. Since $\alpha_l = 0$, we must have $f_1 \circ 
\Psi_{ij}(z) \neq 0$. Hence the point $\Psi_{ij}(z) \in V$ is in the zero set of $f_2$, but not
in the zero set of $f_1$, a contradiction. We conclude that Lojasiewicz's inequality holds and 
we are done.

We next prove Theorem 1.1. Let $f_1(x),...,f_m(x)$ be  
real-analytic functions defined on a neighborhood $V$ of a compact subset $K$ of $\R^n$. Similar
to the proof of Lojasiewicz's inequality, for each $x \in K$ with $\prod_{l=1}^m f_l(x) = 0$, 
let $\phi^x$ 
be a bump function supported in $V$ for which the corollary to the Main Theorem holds for each
$f_l$. If $\prod_{l=1}^m f_l(x) \neq 0$, Let $\phi^x$ be a bump function  
such that $\prod_{l=1}^m f_l \neq 0$ on $supp(\phi^x)$. 
In either case, assume $\phi^x = 1$ on some neighborhood $V_x$ of $x$. By compactness, we
may let $\{V_{x_j}\}$ be finitely many of these sets covering $K$. The set $V' = \cup_jV_{x_j}$
will satisfy the conclusions of the Main Theorem. 

We create a partition of unity based on the $\phi^{x_j}$ by letting $\alpha^j = {\phi^{x_j} \over
\sum_j \phi^{x_j}}$. Thus $\alpha^j \in C^{\infty}(V')$ and $\sum_j \alpha^j = 1$ on $V'$. Hence if
$O$ is any open set with $K \subset O \subset V'$, we have
$$ \int_O \prod_{l=1}^m |f_l(x)|^{-\epsilon_l}\,\,dx = \sum_j \int_O \prod_{l=1}^m |f_l(x)|
^{-\epsilon_l}\alpha^j(x)\,\,dx \eqno (5.2)$$
Decomposing $\phi^{x_j}$ as in the Main Theorem, write $\phi^{x_j} = \sum_i \phi^{x_j}_i$.
Analogously, write $\alpha^j_i = {\phi^{x_j}_i \over \sum_j \phi^{x_j}}$. Thus $\sum_i
\alpha^j_i = \alpha^j$, and we have
$$ \int_O \prod_{l=1}^m |f_l(x)|^{-\epsilon_l}\,\,dx = \sum_{i,j} \int_O \prod_{l=1}^m 
|f_l(x)|^{-\epsilon_l}\alpha^j_i(x)\,dx \eqno (5.3)$$
Let $\Psi_{ij}$ be the composition of the coordinate changes on $\{x: \phi^{x_j}_i > 0\}$ given
by the Main Theorem, but shifted so that $\Psi_{ij}(0) = x_j$ (i.e. instead of  $\Psi_{ij}(0) = 0$).
For the $\phi^{x_j}$ on whose support $\prod_{l=1}^m f_l$ doesn't vanish, one can let $\Psi_{ij}$ be the identity 
map for the purposes of the following arguments. Then there are monomials $m_l(x) = \prod_{k=1}^n x_k^{s_{ijkl}}$, and 
functions $c_l(x)$ not vanishing on $\bar O$ such that
$$ f_l \circ \Psi_{ij}(x) = c_l(x)m_l(x) = c_l(x)\prod_{k=1}^n x_k^{s_{ijkl}} \eqno (5.4)$$
Doing a change of variables in a given term of $(5.3)$, we have
$$\int_O \prod_{l=1}^m |f_l(x)|^{-\epsilon_l}\alpha^j_i(x)dx $$
$$= \int_{D_{ij}} \prod_{l=1}^m |c_l \circ \Psi_{ij}(x)|^{-\epsilon_l}\prod_{k=1}^n x_k^{\sum_{l=1}^m -\epsilon_l 
s_{ijkl}} (\alpha^j_i \circ \Psi_{ij}(x)) |det \,\Psi_{ij}(x)| \,dx \eqno (5.5)$$
Here $|det \,\Psi_{ij}(x)|$ denotes the Jacobian of the coordinate change $\Psi_{ij}$, and $D_{ij}$ is 
an open set whose existence is guaranteed by the Main Theorem such that $\Psi_{ij}$ is an 
isomorphism from $D_{ij}$ to $O \cap \{x: \alpha_i^j(x) > 0\}$. Next, By 3) of the Main Theorem,
there is a ball $B_{ij}$ centered at the origin such that $\Psi_{ij}$ extends to 
$B_{ij}$ with $\Psi_{ij}(0) = x_j$. Shrinking the $B_{ij}$ if necessary, we may assume that
each $B_{ij}$ is the same ball $B$. In addition, since $\Psi_{ij}(0) = x_j \in K \subset O$, we can 
also assume that $B$ is small enough that each $B \subset \Psi_{ij}^{-1}(O)$. Since each 
$g_{ij}^k$ is $k$-to-1 outside a set of measure zero for some $k$, there is some $N$ such that each 
$\Psi_{ij}$ is an $N$ to 1 map from $D_{ij} \cup B$ into $O$ outside a set of measure zero. Consequently, doing a change of
coordinates, there is a function $\gamma_{ij}(x)$ with $1 \leq |\gamma_{ij}(x)| \leq N$ such that
$$\int_O \prod_{l=1}^m |f_l(x)|^{-\epsilon_l}\alpha^j_i(x)\gamma_{ij}(x)\,\,dx =$$
$$\int_{D_{ij} \cup B} \prod_{l=1}^m |c_l \circ \Psi_{ij}(x)|^{-\epsilon_l}\prod_{k=1}^n 
|x_k|^{\sum_{l=1}^m -\epsilon_l s_{ijkl}} (\alpha^j_i \circ \Psi_{ij}(x))|det \,\Psi_{ij}(x)| \,dx$$
As a result we have 
$$ \int_O \prod_{l=1}^m |f_l(x)|^{-\epsilon_l}\alpha^j_i(x) dx \leq 
\int_O \prod_{l=1}^m |f_l(x)|^{-\epsilon_l}\alpha^j_i(x)\gamma_{ij}(x)\,\,dx$$
$$= \int_{D_{ij} \cup B} \prod_{l=1}^m |c_l \circ \Psi_{ij}(x)|^{-\epsilon_l}\prod_{k=1}^n 
|x_k|^{\sum_{l=1}^m -\epsilon_l s_{ijkl}} (\alpha^j_i \circ \Psi_{ij}(x))|det \,\Psi_{ij}(x)| \,dx $$
$$\leq \int_{D_{ij} \cup B} \prod_{l=1}^m |c_l \circ \Psi_{ij}(x)|^{-\epsilon_l}\prod_{k=1}^n 
|x_k|^{\sum_{l=1}^m -\epsilon_l s_{ijkl}}(\alpha^j_i \circ \Psi_{ij}(x) + \chi_B(x))|
det \,\Psi_{ij}(x)| \,dx \eqno (5.6)$$
But because each $|c_l|$ is bounded above and below and each $|\gamma_{ij}(x)|$ is bounded, changing
variables back in the last equation of $(5.6)$ we get
$$\int_{D_{ij} \cup B} \prod_{l=1}^m |c_l \circ \Psi_{ij}(x)|^{-\epsilon_l}\prod_{k=1}^n 
|x_k|^{\sum_{l=1}^m -\epsilon_l s_{ijkl}}(\alpha^j_i \circ \Psi_{ij}(x) + \chi_B(x))|
det \,\Psi_{ij}(x)| \,dx$$
$$\leq  C \int_O \prod_{l=1}^m |f_l(x)|^{-\epsilon_l}\,\,dx \eqno (5.7)$$
Adding $(5.6)$ and $(5.7)$ over all $i,j$ and using $(5.3)$ we have that there is a constant
$C'$  such that
$$\int_O \prod_{l=1}^m |f_l(x)|^{-\epsilon_l}\,\,dx $$
$$\leq \sum_{i,j} \int_{D_{ij} \cup B} \prod_{l=1}^m |c_l \circ \Psi_{ij}(x)|^{-\epsilon_l}\prod_{k=1}^n 
|x_k|^{\sum_{l=1}^m -\epsilon_l s_{ijkl}}(\alpha^j_i \circ \Psi_{ij}(x) + \chi_B(x))|
det \,\Psi_{ij}(x)| \,dx$$
$$ \leq  C' \int_O \prod_{l=1}^m |f_l(x)|^{-\epsilon_l}\,\,dx $$
Since $|c_l \circ \Psi_{ij}(x)|$ is bounded above and below, we conclude that
$$\int_O \prod_{l=1}^m |f_l(x)|^{-\epsilon_l}  < \infty $$
$$\hbox { iff }
\int_{D_{ij} \cup B} \prod_{k=1}^n |x_k|^{\sum_{l=1}^m -\epsilon_l s_{ijkl}} 
(\alpha^j_i \circ \Psi_{ij}(x) + \chi_B(x)) |det \,\Psi_{ij}(x)|dx < \infty \hbox { for all } i, j\eqno (5.8)$$
By the Main Theorem, we can write  $|det \,\Psi_{ij}(x)| = c_{ij}(x) \prod_{k} x_k^{t_{ijk}}$, 
where the 
$|c_{ij}(x)|$ are bounded above and below. Hence   $(5.8)$ implies that
$$\int_O\prod_{l=1}^m |f_l(x)|^{-\epsilon_l}  < \infty $$
$$\hbox { iff } \int_{D_{ij} \cup B} \prod_{k=1}^n |x_k|^{\sum_{l=1}^m -\epsilon_l s_{ijkl} 
+ t_{ijk}} (\alpha^j_i \circ \Psi_{ij}(x) + \chi_B(x)) < \infty \hbox { for all } i, j\eqno (5.9)$$
If $\sum_{l=1}^m -\epsilon_l s_{ijkl} + t_{ijk} > -1$ for each $k$, then each $x_k$ in $(5.9)$ 
appears to a power greater than -1 and the right-hand integral of $(5.8)$ is finite being over a 
bounded domain. On the other
hand, if for some $k$ we have  $\sum_{l=1}^m -\epsilon_l s_{ijkl} + t_{ijk} < -1$, then since 
the domain $D_{ij} \cup B$ contains the ball $B$ centered at the origin on which the
integrand is at least 1, the integral 
in the $x_k$ variable in $(5.8)$ is infinite. Hence the $i, j$ term of the right-hand side of 
$(5.9)$ is finite iff $\sum_{l=1}^m -\epsilon_l s_{ijkl} + t_{ijk} > -1 $ for all $k$. We conclude
that 
$$\int_O\prod_{l=1}^m |f_l(x)|^{-\epsilon_l}  < 
\infty \hbox { iff } \sum_{l=1}^m \epsilon_l s_{ijkl} <  t_{ijk} + 1  \hbox { for all } i, j, 
\hbox { and } k. \eqno (5.10)$$
The equations of $(5.10)$ are independent of $O$, so we have proved Theorem 1.1. It is worth
noting that $(5.10)$ ensures that $\int_O\prod_{l=1}^m |f_l(x)|^{-\epsilon_l}$ is finite whenever 
$\epsilon_1,..., \epsilon_m$ are sufficiently small positive numbers.

{\noindent \bf References:}

\parskip=6pt

\noindent [A] S. Abhyankar, {\it Good points of a hypersurface}, Adv. in Math. {\bf 68} (1988), no. 2,
87-256. 

\noindent [At] M. Atiyah, {\it Resolution of singularities and division of distributions},
Comm. Pure Appl. Math. {\bf 23} (1970), 145-150.

\noindent [BGe] I. N Bernstein and S. I. Gelfand, {\it Meromorphy of the function ${\rm P}^
\lambda$}, Funkcional. Anal. i Prilo\v zen. {\bf 3} (1969), no. 1, 84-85.

\noindent [BM1] E. Bierstone, P. Milman, {\it Semianalytic and subanalytic sets}, Inst. Hautes
Etudes Sci. Publ. Math. {\bf 67} (1988) 5-42.

\noindent [BM2] E. Bierstone, P. Milman, {\it Canonical desingularization in 
characteristic zero by blowing up the maximum strata of a local invariant.} 
Invent. Math. {\bf 128} (1997), no. 2, 207-302.

\noindent [BM3] E. Bierstone, P. Milman, {\it Resolution of singularities in Denjoy-Carleman 
classes.} Selecta Math. (N.S.) {\bf 10} (2004), no. 1, 1-28. 

\noindent [CGiO] V. Cossart, J. Giraud, U Orbanz, {\it Resolution of surface singularities}
(appendix by H. Hironaka). Lecture Notes in Mathematics,  Springer-Verlag, Berlin {\bf 1101},
(1984). 

\noindent [EV1] S. Encinas, O. Villamayor, {\it Good points and constructive resolution of 
singularities.} Acta Math. {\bf 181} (1998), no. 1, 109-158. 

\noindent [EV2] S. Encinas, O. Villamayor, {\it A new proof of desingularization over fields of 
characteristic zero.} Proceedings of the International Conference on Algebraic Geometry and 
Singularities (Spanish) (Sevilla, 2001). Rev. Mat. Iberoamericana {\bf 19} (2003), no. 2, 339-353. 

\noindent [G1] M. Greenblatt, {\it A direct resolution of singularities for
functions of two variables with applications to analysis}, J. Anal. Math. 
{\bf 92} (2004), 233-257.

\noindent [G2] M. Greenblatt, {\it Sharp $L\sp 2$ estimates for one-dimensional 
oscillatory integral operators with $C\sp \infty$ phase.} Amer. J. Math. 
{\bf 127} (2005), no. 3, 659-695. 

\noindent [G3] M. Greenblatt, {\it Newton polygons and local integrability
of negative powers of smooth functions in the plane}, Trans. Amer. Math. Soc. 
{\bf 358} (2006), no. 2, 657-670.

\noindent [G4] M. Greenblatt, {\it Resolution of singularities, asymptotic expansions of
integrals over sublevel sets, and applications}, preprint, to be submitted.

\noindent [G5] M. Greenblatt, {\it Oscillatory integral decay, sublevel set growth,
and the Newton polyhedron}, submitted.

\noindent [H1] H. Hironaka, {\it Resolution of singularities of an algebraic 
variety over a field of characteristic zero I},  Ann. of Math. (2) {\bf 79}
(1964), 109-203;

\noindent [H2] H. Hironaka, {\it Resolution of singularities of an algebraic 
variety over a field of characteristic zero II},  Ann. of Math. (2) {\bf 79}
(1964), 205-326. 

\noindent [K] J. Kollar, {\it Resolution of Singularities - Seattle Lectures}, preprint.

\noindent [L] S. Lojasiewicz, {\it Ensembles semi-analytiques}, Inst. Hautes 
Etudes Sci., Bures-sur-Yvette, 1964.

\noindent [PS] D. H. Phong, E. M. Stein, {\it The Newton polyhedron and
oscillatory integral operators}, Acta Mathematica {\bf 179} (1997), 107-152.

\noindent [PSSt] D. H. Phong, E. M. Stein, J. Sturm, {\it On the growth and 
stability of real-analytic functions}, Amer. J. Math. {\bf 121} (1999), no. 3, 519-554.

\noindent [PSt] D. H. Phong, J. Sturm {\it Algebraic estimates, stability of local zeta 
functions, and uniform estimates for distribution functions}, Ann. of Math. (2) {\bf 152}
(2000), no. 1, 277-329. 

\noindent [R] V. Rychkov, {\it Sharp $L^2$ bounds for oscillatory
integral operators with $C^\infty$ phases}, Math. Zeitschrift, {\bf 236}
(2001) 461-489.

\noindent [S] E. Stein, {\it Harmonic analysis; real-variable methods,
orthogonality, and oscillatory integrals}, Princeton Mathematics Series Vol. 
43, Princeton University Press, Princeton, NJ, 1993.

\noindent [Su] H. Sussman, {\it Real analytic desingularization and subanalytic sets:
an elementary approach}, Trans. Amer. Math. Soc. {\bf 317} (1990), no. 2, 417-461.

\noindent [T] G. Tian, {\it On Calabi's conjecture for complex surfaces with positive 
first Chern class}, Invent. Math. {\bf 101} (1990), no. 1, 101-172.

\noindent [V] A. N. Varchenko, {\it Newton polyhedra and estimates of
oscillatory integrals}, Functional Anal. Appl. {\bf 18} (1976), no. 3, 
175-196.

\noindent [W] J. Wlodarczyk, {\it Simple Hironaka resolution in characteristic zero.}
J. Amer. Math. Soc. {\bf 18} (2005), no. 4, 779-822.

\noindent 244 Mathematics Building \hfill \break
\noindent University at Buffalo \hfill \break
\noindent Buffalo, NY 14260 \hfill \break
\noindent mg62@buffalo.edu
\end